\documentclass[12pt]{article}
\usepackage{fullpage}
\pdfoutput=1
\usepackage{graphicx,dsfont,amsmath,amssymb,authblk} %
\usepackage{color}
\usepackage{epstopdf}
\usepackage{hyperref} 
\usepackage{cite}
 
\DeclareMathOperator{\tr}{tr}
\newcommand{\jmat}{\mathbf{j}}
\newcommand{{\diag}}{\rm{diag}}
\DeclareMathOperator{\sdet}{sdet}
\DeclareMathOperator{\str}{str}
\newcommand{\eins}{\leavevmode\hbox{\small1\kern-3.8pt\normalsize1}}

\begin{document}
\allowdisplaybreaks 

\title{The Correlated Jacobi and the Correlated Cauchy-Lorentz ensembles}
\author{Tim Wirtz$^1$, Daniel Waltner$^1$, Mario Kieburg$^2$, Santosh Kumar$^3$}
\affil{$^1$Fakult\"at f\"ur Physik, Universit\"at Duisburg--Essen, Duisburg, Germany\\
$^2$Fakult\"at f\"ur Physik, Universit\"at Bielefeld, Bielefeld, Germany\\
$^3$Department of Physics, Shiv Nadar University, Gautam Buddha Nagar, Uttar Pradesh - 201314, India}

\date{}
\maketitle
\begin{abstract}
We calculate the $k$-point generating function of the correlated Jacobi ensemble using supersymmetric methods. We use the result for complex matrices for $k=1$ to derive a closed-form expression 
for eigenvalue density. For real matrices we obtain the density in terms of a twofold integral that we evaluate numerically. For both expressions we find agreement when comparing with Monte Carlo simulations. Relations between these quantities for the Jacobi and the Cauchy-Lorentz ensemble are derived.
\end{abstract}

\section{Introduction}

The Jacobi ensemble, like the Wishart ensemble, has its roots in the field of multivariate statistics. To quantify the empirically estimated canonical correlation coefficients between two sets of time series, they are compared to a null hypothesis, i.e. to a Gaussian distribution. Assuming for both sets Gaussian statistics with a non-trivial correlation structure, the null hypothesis becomes the Jacobi model~\cite{anderson,muirhead,Johns}. Besides the canonical correlation analysis it applies also to other aspects of high-dimensional inference such as analysis of variances, regression analysis and the test of equality of covariance matrices~\cite{anderson,MardiaKentBibby1979,Johns,muirhead,JohnsonWichern2013}. Consequently, they are also referred to as MANOVA ({\bf M}ultivariate {\bf AN}alysis {\bf O}f {\bf VA}riance) ensembles. 

Along with the Gaussian (Wigner) and the Laguerre (uncorrelated-Wishart) ensembles, the Jacobi ensemble constitutes the family of classical random matrix ensembles. Just as the eigenvalue statistics for Wigner and Wishart ensembles respectively involve Hermite polynomials and associated Laguerre polynomials, the eigenvalue statistics of Jacobi ensembles involve Jacobi polynomials. Together they complete the random matrix ensemble picture in connection with the theory of classical orthogonal, and skew-orthogonal polynomials~\cite{AFM2009,GP2002,KP2011,ForresterLogGas}. 
Closely related to these classical ensembles is the less known Cauchy-Lorentz ensemble~\cite{KaymakKieburgGuhr,Kie15}. It exhibits a Levy tail and therefore finds important applications in the spectral statistics of covariance matrices in finance~\cite{BJNPZ}.

Besides their crucial role in multivariate statistics and the intimate connection with the theory of classical polynomials, Jacobi ensembles find interesting applications in the fields of quantum transport and optical fibre communication. In the context of quantum transport, they describe the statistics of transmission (or reflection) eigenvalues for disordered mesoscopic systems with {\it ideal leads}, see \cite{Beenakker1997} and references therein. The knowledge of these eigenvalues, in turn, gives access to important observables such as Landauer conductance, shot-noise power and Wigner delay time~\cite{Forrester2006,SS2006,SSW2008,KP2010,VMB2010,MS2011}. This connection of Jacobi ensembles with the transmission eigenvalues stems from their relationship with the scattering matrices which are modelled using Dyson's circular ensembles~\cite{Forrester2006,Beenakker1997}. More recently Jacobi ensembles have been used in the ergodic capacity analysis of multiple-input-multiple-output optical fibre 
communication~\cite{DarFederShtaif, DFS2013,KMV2014}. In this case the Jacobi structure emerges from the channel matrix which happens to be a block of a bigger transfer matrix which is unitary in nature. 

 In all these applications listed above one can ask about intrinsic correlations of the channels or time series which result in correlated random matrix ensembles. Thus one has to relax the condition of the empirical correlation matrix being proportional to the identity matrix. This correlation can be attributed to different reasons depending on the context. For instance, in time-series analysis data are in general mutually correlated, e.g. see~\cite{Toole}. In the context of multiple antenna communication, this correlation arises because of spatial correlation between closely spaced antennas~\cite{ATLV} as it is the case in any cellphone. To handle these situations matrix models with correlations between different time series~\cite{ATLV,Rech12,For06,WKG14,WirtzKieburgGuhrPRL}, even with double correlations~\cite{SimMou,WaltnerWirtzGuhr2014}, have been introduced; see also Refs.~\cite{anderson,muirhead} and references therein. In the Wishart case the level density of the model correlation matrix was computed by 
various techniques. We pursue the ideas of Refs.~\cite{Rech12,WaltnerWirtzGuhr2014} where supersymmetry techniques were employed.

In the case of the Jacobi and the Cauchy-Lorentz ensemble not much progress has been made in a similar direction as  for the Gaussian because of their non-Gaussian form. Nonetheless both models can be traced back to a combination of two random matrices both drawn from Wishart ensembles, as explained in detail in sections~\ref{sec-Jac} and~\ref{sec:JCL}. To each of the two Wishart matrices one can associate an empirical covariance matrix.  However, we show in section~\ref{sec-Jac} that effectively only one empirical correlation matrix is involved. 

In section~\ref{sec-Jac} we briefly review the Jacobi ensemble before showing the relation to the Cauchy-Lorentz ensemble in section~\ref{sec:cje:lorentzianmodel}. The latter is more convenient to apply the projection formula~\cite{KaymakKieburgGuhr} and the supersymmetry approach. Nevertheless we present in appendix~\ref{sec:cje:twosupermatrixmodel} an alternative approach which is a direct way from the integral over the Jacobi ensemble to an integral over superspace. The resulting supersymmetric integral is more involved than the one resulting from the first approach because it involves two supermatrices instead of only one. However the second approach exhibits an intrinsic symmetry of the model which is not obvious in the more compact result presented in section~\ref{sec:cje:lorentzianmodel}. In section~\ref{sec:cje:complexcase} we calculate the whole eigenvalue statistics of the correlated complex Jacobi and the correlated complex Cauchy-Lorentz ensemble. This is possible due to an underlying 
determinantal point process. Such an integrable structure is not available for the real case discussed in section~\ref{sec:cje:realcase}. For this case we need the supersymmetric representation and compute the level density with the help of the generalized Hubbard-Stratonovich transformation. Some details of this calculation are presented in appendix~\ref{app:expl}. In section \ref{sec:saddlepoint}, we discuss the level density for the real and the complex ensemble in the limit of large matrix sizes in a unified way. Our results are summarized in section \ref{sec:conclusion}.

\section{Correlated Jacobi Ensemble}\label{sec-Jac}

This section is devoted to introduce the theoretical background.  The correlated Jacobi ensemble is a two-matrix model for a Hermitian matrix $ H$ which is composed of two independently distributed correlated Wishart matrices $FF^\dagger$ and $BB^\dagger$, i.e.
\begin{align}
 \label{eq:cje:defH}
 H &= \frac{FF^\dagger - BB^\dagger}{FF^\dagger + BB^\dagger}.
\end{align}
The rectangular matrices $F$ and $B$ have the same number of rows $p$ (the number of time series) but may have a different number of columns $n_1$ and $n_2$ (the number of time steps), respectively. We assume that $n_1\geq p$ and $n_2\geq p$. In a realistic situation for time series analysis both dimensions $n_1$ and $n_2$ are indeed larger than $p$ since the number of time steps is usually larger than the number of time series. Moreover the matrix entries $F_{ij}$  are distributed as Gaussians with variances $\left[C_F\right]_{ji}$ such that the distribution reads
\begin{align}
\label{dist}
P(F|C_F)=\frac{\exp\left(-\tr FF^\dagger C_F^{-1}/\gamma\right)}{\left(\gamma\pi\right)^{n_1p/\gamma}\det^{n_1p/\gamma_1}C_F},
\end{align}
and likewise for $B$, where we introduce $\gamma=1$ for $\beta=2$ and $\gamma=2$ for $\beta=1$. The parameter $\beta$ is the Dyson index meaning that $\beta=1$ corresponds to real matrices whereas $\beta=2$ denotes complex matrices. The case of quaternion matrices ($\beta=4$) is not considered here but can be in general worked out in a similar way.

The empirical correlation matrices $C_F$ and $C_B$ are fixed as in the discussions of \cite{ATLV,SimMou,For06,Rech12,WaltnerWirtzGuhr2014,WKG14,WirtzKieburgGuhrPRL}. Because of the distribution~\eqref{dist} the random matrices $F$ and $B$ have upon average the same covariances as the sample,
\begin{equation}
\frac{1}{n_1}\left\langle FF^\dagger\right\rangle=C_F,\quad \frac{1}{n_2}\left\langle BB^\dagger\right\rangle=C_B
\end{equation}
The angular brackets denote the average over the random matrices $F$ and $B$.

We note that instead of the matrix model (\ref{eq:cje:defH}) one can also consider the matrix model $FF^\dagger(FF^\dagger+BB^\dagger)^{-1}$. In that case the eigenvalues lie in the interval $[0,1]$. These two matrix models are related by a scaling and a shift transformation. We work here with the matrix model (\ref{eq:cje:defH}) and thus relate directly to the Jacobi weight. Because of the involvement of two Wishart matrices in the Jacobi model, it is also known as the double Wishart model.

The spectral statistics of $H$ can be studied via partition functions which are averages over products and ratios of characteristic polynomials, see e.g.\ \cite{GuhrHabil,Haake}, 
\begin{align}
 \label{eq:cje:reminderGeneratingFunction}
 Z_{p,\beta}^{k_1|k_2}(\kappa) &= \int d[F,B] \frac{\prod_{a=1}^{k_2} \det\left(H - \kappa_{a2}\eins_{ p}\right)}{\prod_{b=1}^{k_1} \det\left(H- \kappa_{b1}
 \eins_{p}\right)}
 P(F|C_F)P(B|C_B).
\end{align}
Those partition functions are related to the $k$-point correlation functions as
\begin{eqnarray}
 R_{p,\beta}^{k}(x) &=&\left\langle\prod_{j=1}^k\left(\frac{1}{p}\tr\delta(H-x_j\eins_p)\right)\right\rangle \label{eq:k-point-self}\\
 & =&\left.\lim_{\epsilon\to0}\sum_{L_j=\pm1}\left(\prod_{j=1}^k\frac{L_j}{2\pi\imath p}\frac{\partial}{\partial x'_j}\right) Z_{p,\beta}^{k|k}(x_1+\imath L_1\epsilon,\ldots,x_k+\imath L_k\epsilon,x'_1,\ldots,x'_k)\right|_{x'=x}
\nonumber
\end{eqnarray}
with $x_1,\ldots,x_k\in\mathbb{R}$. We underline that this definition of the $k$-point correlation function comprises self-energy terms, in particular terms proportional to $\delta(x_i-x_j)$. Despite those terms the definition~\eqref{eq:k-point-self} is very helpful when the joint probability density of the eigenvalues is not explicitly known as it is in the case of the correlated real ensembles discussed in section~\ref{sec:cje:realcase}.

In the case when  the joint probability density is known one can apply another definition of the $k$-point correlation function, e.g. see \cite{Mehta},
\begin{eqnarray}
\widehat{R}_{p,\beta}^{k}(x)&=&\left\langle\prod_{j=1}^k\delta(E_j-x_j)\right\rangle=\frac{(p-k)!p^k}{p!}R_{p,\beta}^{k}(x)-\textrm{ self-energy\ terms},
 \label{eq:k-point}
\end{eqnarray}
where $E_1,\ldots,E_k$ are the eigenvalues of $H$. This definition is employed to the correlated complex ensembles studied in section~\ref{sec:cje:complexcase}. Both definitions are normalized in such a way that $\int dx_k R_{p,\beta}^{k}(x)=R_{p,\beta}^{k-1}(x)$ and $\int dx_k \widehat{R}_{p,\beta}^{k}(x)=\widehat{R}_{p,\beta}^{k-1}(x)$. The eigenvalue density  of $\mathcal{H}$ is given by the case $k=1$, i.e.
\begin{equation}\label{density-def}
S_\beta(x):=\left\langle\frac{1}{p}\tr\delta(H-x\eins_p)\right\rangle=R_{p,\beta}^{1}(x)=\widehat{R}_{p,\beta}^{1}(x).
\end{equation}
 It is computed in sections~\ref{sec:cje:complexcase} and \ref{sec:cje:realcase}.

The partition function~\eqref{eq:cje:reminderGeneratingFunction} is normalized as $\lim_{\kappa\to\infty}Z_{p,\beta}^{k_1|k_2}(\kappa){\det}^p(-\kappa_1)/{\det}^p(-\kappa_2)\rightarrow 1$. Note that in the case $C_F=C_B$, the partition function and, hence, the spectral statistics of $H$ become independent of the empirical covariance matrices such that  the joint eigenvalue distribution of $H$ is given by the uncorrelated Jacobi ensemble \cite{muirhead,MardiaKentBibby1979,ForresterLogGas}. Indeed we can rescale the matrices $F\to\sqrt{C_F}F$ and $B\to\sqrt{C_F}B$ such that the partition function only depends on $C_{\text{eff}}=C_F^{-1/2}C_BC_F^{-1/2}$, i.e.
\begin{align}
 \label{eq:cje:GeneratingFunctionCeff}
 Z_{p,\beta}^{k_1|k_2}(\kappa) &= \int d[F,B] \frac{\prod_{a=1}^{k_2} \det\left(H - \kappa_{a2}\eins_{p}\right)}{\prod_{b=1}^{k_1} \det\left(H- \kappa_{b1}\eins_{p}\right)}
 P(F|\eins_{p})P(B|C_{\text{eff}}).
 \end{align}
 Thus in the case $C_F\neq C_B$, the eigenvalue statistics become non-trivial, even in the complex case ($\beta=2$), see section~\ref{sec:cje:complexcase}. Thus, it is reasonable to apply other methods than the standard Jack or Zonal polynomial approach \cite{muirhead,Dumitriuetal,dubbsetal}.

\section{Correlated Cauchy-Lorentz Ensemble}\label{sec:cje:lorentzianmodel}

The correlated Jacobi ensemble is also related to the correlated Cauchy-Lorentz ensemble. This relation is established in subsection~\ref{sec:JCL}. With the help of this relation we express the partition function~\eqref{eq:cje:GeneratingFunctionCeff} in terms of integrals over supermatrices in subsection~\ref{sec:proj}. The advantage of the supersymmetric integrals is the drastic reduction of integration variables. The parameters $n_1$, $n_2$, and $p$ only appear as external parameters in those expressions. This fact allows an asymptotic study when these parameters are large  which is considered in section~\ref{sec:saddlepoint}.

\subsection{Relation between Jacobi and Cauchy-Lorentz}\label{sec:JCL}

To see the relation between the Jacobi and the Cauchy-Lorentz ensembles we take advantage of the fact that $FF^\dagger$ is generically invertible because $n_1>p$. Then we can rewrite $F=\widehat{F}\Pi$ as a product of a $p\times p$ square matrix $\widehat{F}$ which can be either real ($\beta=1$) or complex ($\beta=2$) and a random projection $\Pi\in{\rm O}(n_1)/[{\rm O}(p)\times{\rm O}(n_1-p)]$  for $\beta=1$ and $\Pi\in{\rm U}(n_1)/[{\rm U}(p)\times{\rm U}(n_1-p)]$  for $\beta=2$. The measure transforms as $d[F]={\det}^{(n_1-p)/\gamma}\widehat{F}\widehat{F}^\dagger d[\widehat{F}]d\mu(\Pi)$ where $d\mu(\Pi)$ is the Haar measure induced   from those on the groups ${\rm O}(n_1)$ or ${\rm U}(n_1)$, respectively. The projection drops out and the integral over it yields a constant. The same procedure can be applied for the matrix $B$ yielding a $p\times p$ matrix $\widehat{B}$ with the weight $ \exp[-\tr C_{\rm eff}^{-1}\widehat{B}\widehat{B}^\dagger/\gamma]{\det}^{(n_2-p)/\gamma}\widehat{B}\widehat{B}^\dagger d[\
widehat{B}]$

Since $\widehat{F}$ is generically invertible we can rewrite any average of an observable $\mathcal{O}$ of $H$ in the following way
\begin{eqnarray}
\left\langle\mathcal{O}\left(\frac{FF^\dagger-BB^\dagger}{FF^\dagger+BB^\dagger}\right)\right\rangle=\left\langle\mathcal{O}\left(\frac{\widehat{F}\widehat{F}^\dagger -\widehat{B}\widehat{B}^\dagger}{\widehat{F}\widehat{F}^\dagger+\widehat{B}\widehat{B}^\dagger}\right)\right\rangle=\left\langle\mathcal{O}\left(\frac{\eins_p -\widehat{F}^{-1}\widehat{B}\widehat{B}^\dagger\widehat{F}^{-1\,\dagger}}{\eins_p+\widehat{F}^{-1}\widehat{B}\widehat{B}^\dagger\widehat{F}^{-1\,\dagger}}\right)\right\rangle.\label{average-1}
\end{eqnarray}
Note that the observable $\mathcal{O}$ has to be invariant under the group action of ${\rm O}(p)$ and ${\rm U}(p)$, respectively. Due to  the invariance of the probability weight $P(F|\eins_{p})$ the average~\eqref{average-1} over $\widehat{F}$, only, depends on the singular values of $\widehat{B}$ such that it is also true
\begin{eqnarray}
\left\langle\mathcal{O}\left(\frac{FF^\dagger-BB^\dagger}{FF^\dagger+BB^\dagger}\right)\right\rangle=\left\langle\mathcal{O}\left(\frac{\eins_p -\widehat{F}^{-1}\widehat{B}^\dagger\widehat{B}\widehat{F}^{-1\,\dagger}}{\eins_p+\widehat{F}^{-1}\widehat{B}^\dagger\widehat{B}\widehat{F}^{-1\,\dagger}}\right)\right\rangle.\label{average-2}
\end{eqnarray}
This identity is reminiscent of the weak-commutation relation proven in~\cite{IpsKie}. When rescaling $\widehat{B}\to\widehat{B}\widehat{F}^\dagger$ the integral over $\widehat{F}$ becomes a deformed Gaussian and reads
\begin{eqnarray}
&&\int d[\widehat F]{\det}^{(n_1+n_2-p)/\gamma}\widehat F\widehat F^\dagger \exp\left(-\frac{1}{\gamma}\tr \widehat F\left[ \widehat{B}^\dagger C_{\text{eff}}^{-1} \widehat{B} +\eins_p\right]\widehat{F}^\dagger\right)\nonumber\\
&\propto&{\det}^{-(n_1+n_2)/\gamma}\left(\widehat{B}^\dagger C_{\text{eff}}^{-1} \widehat{B} +\eins_p\right)~.\label{eq:cje:afterBtoBF}
\end{eqnarray}
Thus the spectral statistics of $H$ is equal to the statistics  of
\begin{eqnarray}\label{Hprime}
 H'=\frac{\eins_p-B'B'^\dagger}{\eins_p+B'B'^\dagger},
\end{eqnarray}
where the $p\times n_2$ matrix $B'$ is drawn from the correlated Cauchy-Lorentz distribution
\begin{eqnarray}
 P_{n_1+n_2}^{\rm CL}(B'|C_{\rm eff})&=&\pi^{n_2p/\gamma} \left(\prod_{j=0}^{n-1}\frac{\Gamma[(n_1-p+j+1)/\gamma]}{\Gamma[(n_1+n_2-p+j+1)/\gamma]}\right){\det}^{n_1/\gamma}C_{\text{eff}}\nonumber\\
 &&\times{\det}^{-(n_1+n_2)/\gamma}\left(B'B'^\dagger + C_{\text{eff}}\right).\label{defCL}
\end{eqnarray}
Indeed the factorization $B'=\widehat{B}\Pi'$ in a square matrix $\widehat{B}$ and a projection $\Pi'\in{\rm O}(n_2)/[{\rm O}(p)\times{\rm O}(n_2-p)]$  for $\beta=1$ and $\Pi'\in{\rm U}(n_2)/[{\rm U}(p)\times{\rm U}(n_2-p)]$  for $\beta=2$ is still possible. Then we would have the additional term ${\det}^{(n_2-p)/\gamma}\widehat{B}\widehat{B}^\dagger$ in the weight~\eqref{defCL}.

Note that the eigenvalue statistics of $H'$ is completely determined by the eigenvalue statistic of $B'$. This means when we calculate the $k$-point correlation function of $B'B'^\dagger$ then the $k$-point correlation function of $H'$ is given by the substitution $b=(1-x)/(1+x)$ where $b$ is an eigenvalue of $B'B'^\dagger$ and $x$ is an eigenvalue of $H'$. Therefore we consider the partition function
\begin{align}
  \label{eq:cje:lorentziangenfuncII}
\begin{split}
{Z'}_{p,\beta}^{k_1|k_2}(\kappa) &=  \int d[B'] \frac{\prod_{a=1}^{k_2} \det\left(B'B'^\dagger - \kappa_{a2}\eins_{p}\right)}{\prod_{b=1}^{k_1} \det\left(B'B'^\dagger - \kappa_{b1}\eins_{p}\right)}
 P_{n_1+n_2}^{\rm CL}(B'|C_{\rm eff}).
\end{split}
\end{align}
From this point on everything works analogously to the correlated Wishart ensemble studied in \cite{Rech12,WirtzKieburgGuhrPRL}. When plugging Eq.~\eqref{Hprime} into Eq.~\eqref{eq:cje:GeneratingFunctionCeff} we obtain an explicit relation of ${Z'}_{p,\beta}^{k_1|k_2}(\kappa) $ to the partition function of the Jacobi ensemble which is
\begin{eqnarray}\label{relation}
{Z}_{p,\beta}^{k_1|k_2}(\kappa)=(-1)^{(k_1-k_2)p}\frac{\prod_{a=1}^{k_2} (1+\kappa_{b2})^p}{\prod_{b=1}^{k_1} (1+\kappa_{b1})^p}{Z'}_{p,\beta}^{k_1+k_2|k_2+k_1}\left(\frac{1-\kappa_1}{1+\kappa_1},-\eins_{k_2};\frac{1-\kappa_2}{1+\kappa_2},-\eins_{k_1}\right).
\end{eqnarray}
Especially the $k$-point correlation function remains effectively unaffected since in this case we have $k_1=k_2=k$ and the additional characteristic polynomials cancel. The $k$-point correlation function with self-energy terms is then
\begin{eqnarray}
 {R'}_{p,\beta}^{k}(b)\prod_{j=1}^kdb_j &=&\left\langle\prod_{j=1}^k\left(\frac{1}{p}\tr\delta(B'B'^\dagger-b_j\eins_p)db_j\right)\right\rangle \label{eq:k-point-CL-self}\\
 & =&{R'}_{p,\beta}^{k}\left(\frac{1-x}{1+x}\right)\prod_{j=1}^k\frac{2dx_j}{(1+x_j)^2}=R_{p,\beta}^{k}\left(x\right)\prod_{j=1}^kdx_j\nonumber
\end{eqnarray}
and similar for the definition~\eqref{eq:k-point}. We have written the differentials to underline the transformation properties under changes of coordinates. In the particular case of the level density we have
\begin{equation}\label{density-CL-def}
 S'_\beta(b)db=\left\langle\frac{1}{p}\tr\delta(B'B'^\dagger-b\eins_p)db\right\rangle=S'_\beta\left(\frac{1-x}{1+x}\right)\frac{2dx}{(1+x)^2}=S_\beta(x)dx.
\end{equation}
Establishing this relation between the correlated Jacobi model and the correlated Cauchy-Lorentz ensemble is the first main result of our work.

\subsection{Projection Formula}\label{sec:proj}

Before we come to supersymmetry let us refer the reader who is not familiar to superanalysis and superalgebra to the textbook by Berezin~\cite{Berezin}. A general introduction to the supersymmetry method in random matrix theory is given in \cite{Guhrbook,Efebook,GuhrHabil} and references therein.

To apply the projection formula introduced for chiral ensembles in Ref.~\cite{KaymakKieburgGuhr} we rescale the matrix $B'\to \sqrt{C_{\rm eff}}B'$. Then the partition function~\eqref{eq:cje:lorentziangenfuncII} is
\begin{align}
  \label{eq:cje:lorentziangenfuncII-b}
\begin{split}
{Z'}_{p,\beta}^{k_1|k_2}(\kappa) &=  {\det}^{k_2-k_1}C_{\rm eff}\int d[B'] \sdet^{-1}(B'B'^\dagger\otimes\eins_{k_1|k_2}-C_{\rm eff}^{-1}\otimes\kappa)
 P^{\rm CL}(B'|\eins_p).
\end{split}
\end{align}
We use the short-hand notation of the superdeterminant which only encodes the products and ratios of determinants. Moreover we understand $\kappa={\diag}(\kappa_{11},\ldots,\kappa_{k_11},\kappa_{12},\ldots,\kappa_{k_22})$ as a $(k|k)\times(k|k)$ diagonal supermatrix. The representation~\eqref{eq:cje:lorentziangenfuncII-b} directly reflects the duality between the ordinary matrix space and the supermatrix space.

Let $L={\rm sign}\,{\rm Im}\,\kappa$ be the sign of the imaginary increment in the source. Then we can apply the projection formula~\cite{KaymakKieburgGuhr,Kie15}  to the partition function~\eqref{eq:cje:lorentziangenfuncII-b} which reads
\begin{eqnarray}
{Z'}_{p,\beta}^{k_1|k_2}(\kappa)&=& K_1\,{\det}^{k_2-k_1}C_{\rm eff}\int d[\chi] \sdet^{-1/\gamma}(\eins_p\otimes\chi\chi^\dagger-C_{\rm eff}^{-1}\otimes(\kappa\otimes\eins_\gamma))\label{eq:cje:chifinal}\\
&&\times\sdet^{-\mu/\gamma}(\chi\chi^\dagger+\eins_{\gamma k_1|\gamma k_2})\nonumber\\
&=&K_2\,{\det}^{k_2-k_1}C_{\rm eff}\int d\mu(U)\sdet^{n_2/\gamma}U \sdet^{-1/\gamma}(\eins_p\otimes U-C_{\rm eff}^{-1}\otimes(L\kappa\otimes\eins_\gamma))\nonumber\\
&&\times
\sdet^{-\mu/\gamma}((L\otimes\eins_\gamma)U+\eins_{\gamma k_1|\gamma k_2})\nonumber\\
&=&K_3\,{\det}^{k_2-k_1}C_{\rm eff}\int d[\sigma]I_{n_2}(\sigma)\sdet^{-1}\left( \eins_p\otimes\sigma-C_{\rm eff}^{-1}\otimes(\kappa\otimes\eins_\gamma)\right)\nonumber\\
&&\times\sdet^{-\mu/\gamma}\left(\sigma + \eins_{\gamma k_1|\gamma k_2}\right)
\nonumber
\end{eqnarray}
with $\mu=n_1+n_2-p+k_1-k_2$. The normalization constants are
\begin{eqnarray}
 K_1^{-1}&=&\int d[\chi]
\sdet^{-\mu/\gamma}(\chi\chi^\dagger+\eins_{\gamma k_1|\gamma k_2}),\nonumber\\
K_2^{-1}&=&\int d\mu(U) \sdet^{-\mu/\gamma}((L\otimes\eins_\gamma)U+\eins_{\gamma k_1|\gamma k_2})\sdet^{n_2/\gamma}U,\\
K_3^{-1}&=&\int d[\sigma]I_{n_2}(\sigma)\sdet^{-\mu/\gamma}\left(\sigma + \eins_{\gamma k_1|\gamma k_2}\right).\nonumber
\end{eqnarray}
The rectangular supermatrix $\chi$ has dimension $(\gamma k_1|\gamma k_2)\times n_2$ and it satisfies the symmetry
\begin{equation}
\chi^*={\diag}(\eins_{\gamma k_1},\tau_2\otimes\eins_{k_2})\chi,\ (\chi^\dagger)^*=\chi^\dagger{\diag}(\eins_{\gamma k_1},-\tau_2\otimes\eins_{k_2})
\end{equation}
in the real case ($\beta=1$) where $\tau_2$ is the second Pauli matrix, i.e. $\chi$ consists of a $(\gamma k_1)\times n_2$ real or complex matrix depending on $\beta$ and a $(\gamma k_2)\times n_2$ matrix comprising independent Grassmann variables.

In the second equality of Eq.~\eqref{eq:cje:chifinal} we have used the superbosonization formula~\cite{KSG,LSZ,Sommers}. We use this representation in sections~\ref{sec:cje:complexcase} and~\ref{sec:saddlepoint}. In the present situation the matrix
\begin{equation}
 U=\left[\begin{array}{cc}
  U_{\rm BB} & U_{\rm BF} \\ U_{\rm FB} & U_{\rm FF}
 \end{array}\right]
\end{equation}
satisfies the following symmetries:
\begin{itemize}
 \item $(L\otimes\eins_\gamma)U_{\rm BB}=[(L\otimes\eins_\gamma)U_{\rm BB}]^\dagger$ positive definite  for both $\beta=1,2$ and $U_{\rm BB}=U_{\rm BB}^*$ only for $\beta=1$,
 \item $U_{\rm FF}^\dagger=U_{\rm FF}^{-1}$ unitary for both $\beta=1,2$ and $U_{\rm FF}=\tau_2 U_{\rm FF}^T\tau_2$ self-dual only for $\beta=1$,
 \item $U_{\rm BF}=U_{\rm FB}^\dagger$ for both $\beta=1,2$ and $U_{\rm FB}^T=U_{\rm BF}\tau_2$ only for $\beta=1$.
\end{itemize}
The set of supermatrices $U$ is a particular case of the co-set ${\rm Gl}(k|k)/{\rm U}(k-r,r|k)$ for $\beta=2$ and ${\rm UOSP}(2k-2r,2r|2k)/{\rm U}(2k-2r,2r|2k)$ for $\beta=1$ ($r$ is the number of minus signs in $L$), see~\cite{Zirn}.  Note that we have chosen the non-compact group symmetries due to the non-trivial signs $L$.

Alternative to the superbosonization formula one can also choose the generalized  Hubbard-Stratonovich transformation~\cite{IngSie,KGG09,KSG}, see third equality of Eq.~\eqref{eq:cje:chifinal}, which we employ in section~\ref{sec:cje:realcase}. The superbosonization formula and the generalized Hubbard-Stratonovich transformation are equivalent~\cite{KSG}. The first can be understood as the contour representation of the latter which is some kind of a high dimensional residue theorem. The distribution $I_{n_2}(\sigma)$ in the generalized Hubbard-Stratonovich transformation is the supersymmetric Ingham-Siegel integral~\cite{IngSie,KGG09},
\begin{align}
  I_{n_2}(\sigma) = \int d[\rho]\sdet^{-n_2/\gamma}(\rho-\imath\epsilon\eins_{\gamma k_1|\gamma k_2})\exp\left[\imath \str (\rho-\imath\epsilon\eins_{\gamma k_1|\gamma k_2})\sigma\right]\label{eq:gHS:InghamSiegel}
\end{align}
with $\epsilon>0$. It encodes derivatives of Dirac delta functions in the fermion-fermion block, $\sigma_{\rm FF}$, and the positivity condition in the boson-boson block, $\sigma_{\rm BB}$. Therefore the supermatrices $\sigma$ and $\rho$ are elements of another realizations of the cosets ${\rm Gl}(k|k)/{\rm U}(k|k)$ for $\beta=2$ and ${\rm UOSP}(2k|2k)/{\rm U}(2k|2k)$ for $\beta=1$. They are of the form
\begin{equation}
 \sigma=\left[\begin{array}{cc}
  \sigma_{\rm BB} & \sigma_{\rm FB}^\dagger \\ \sigma_{\rm FB} & \imath\sigma_{\rm FF}
 \end{array}\right]\ {\rm and}\ \rho=\left[\begin{array}{cc}
  \rho_{\rm BB} & \rho_{\rm FB}^\dagger \\ \rho_{\rm FB} & \imath\rho_{\rm FF}
 \end{array}\right],
\end{equation}
where $\sigma_{\rm FB}$ and $\rho_{\rm FB}$ comprise independent Grassmann variables, only, and satisfy $\sigma_{\rm FB}^*= \imath\tau_2\sigma_{\rm FB}$ and $\rho_{\rm FB}^*= \imath\tau_2\rho_{\rm FB}$ for $\beta=1$. The blocks $\sigma_{\rm BB}$, $\sigma_{\rm FF}$, $\rho_{\rm BB}$, and $\rho_{\rm FF}$ are Hermitian. Additionally the submatrices $\sigma_{\rm BB}$ and $\rho_{\rm BB}$ are real symmetric while $\sigma_{\rm FF}$ and $\rho_{\rm FF}$ are self-dual for $\beta=1$.

\section{Eigenvalue Spectrum for $\beta=2$}\label{sec:cje:complexcase}

In the case of complex matrices ($\beta=2$) we first calculate the joint probability density before applying supersymmetry. This simplifies the whole calculation a lot. In particular we show in subsection~\ref{sec:det} that the joint probability density follows a determinantal point process. The corresponding kernel is expressed in terms of a supermatrix integral in subsection~\ref{sec:kernel}.

\subsection{Determinantal Point Process}\label{sec:det}

 We start from the weight~\eqref{defCL}. Let $b={{\diag}}(b_1,\ldots,b_p)>0$ be the eigenvalues of $B'B'^\dagger=UbU^\dagger$ with  $U\in{\rm U}(p)$. The computation of the joint probability density of $b$ is possible because the integral
\begin{eqnarray}\label{group-int}
 &&\int_{{\rm U}(p)}d\mu(U){\det}^{-n_1-n_2}\left(UbU^\dagger + C_{\text{eff}}\right)\\
 &=&(-1)^{p(p-1)/2}\left(\prod\limits_{j=0}^{p-1}\frac{j!(n_1+n_2-p-j-1)!}{(n_1+n_2-p-1)!}\right)\frac{\det[(b_i+\Lambda_j)^{p-n_1-n_2-1}]_{1\leq i,j\leq p}}{\Delta_p(b)\Delta_p(\Lambda)}\nonumber
\end{eqnarray}
is well-known~\cite{GroRic,HarOrl}. The matrix $\Lambda={\diag}(\Lambda_1,\ldots,\Lambda_p)>0$ are the eigenvalues of the empirical matrix  $ C_{\text{eff}}$. The measure $d\mu(U)$ is the normalized Haar measure on the group ${\rm U}(p)$ and $\Delta_p(b)=\prod_{1\leq i<j\leq p}(b_j-b_i)=\det[b_i^{j-1}]_{1\leq i,j\leq p}$ is the Vandermonde determinant. Thus the joint probability density is
\begin{equation}\label{jpd}
p_{p,2}^{\rm jpd}(b)=K_{\rm jpd}\,{\det}^{n_1}\Lambda\frac{\Delta_p(b)}{\Delta_p(\Lambda)}\det[b_i^{n_2-p}(b_i+\Lambda_j)^{p-n_1-n_2-1}]_{1\leq i,j\leq p},\ b>0,
\end{equation}
with the normalization constant
\begin{equation}\label{jpd-constant}
K_{\rm jpd}=\frac{1}{p!}\prod\limits_{j=0}^{p-1}\frac{(n_1+n_2-p)!}{(n_2-p+j)!(n_1-j+1)!}.
\end{equation}
The normalization can be readily checked when applying Andr\'eief's integration theorem~\cite{Andr}.

The $k$-point correlation function of the correlated Lorentz-ensemble without self-energy terms is
\begin{equation}\label{k-point-def-prime}
{R'}_{p,2}^k(b_1,\ldots,b_k)= \prod_{j=k+1}^p\int_0^\infty db_j\, p_{p,2}^{\rm jpd}(b).
\end{equation}
We employ a modification of Andr\'eief's integration theorem derived in~\cite{KieGuh} and find a determinantal point process
\begin{eqnarray}
{R'}_{p,2}^k(b_1,\ldots,b_k)&=&\frac{(-1)^k(p-k)! K_{\rm jpd}}{\Delta_p(\Lambda)}\det\left[\begin{array}{cc} 0 & b_m^{n_2-p}(b_m+\Lambda_j)^{p-n_1-n_2-1}\Lambda_j^{n_1} \\ b_l^{i-1} & \displaystyle \frac{(n_2-p+i-1)!(n_1-i)!}{(n_1+n_2-p)!}\Lambda_j^{i-1} \end{array}\right]_{\substack{1\leq l,m\leq k\\1\leq i,j\leq p}}\nonumber\\
&=&\frac{(p-k)!p^k}{p!}\det\left[K(b_l,b_m;\Lambda)\right]_{1\leq l,m\leq k}.\label{k-point}
\end{eqnarray}
We emphasize that the case $k=p$ is the joint probability density. The kernel is
\begin{eqnarray}
 K(b_l,b_m;\Lambda)&=&\frac{1}{p}\sum_{i,j=1}^p (-1)^{p-i}\frac{(n_1+n_2-p)!}{(n_2-p+i-1)!(n_1-i)!} \frac{e_{p-i}^{p-1}(\Lambda_{\neq j})}{\det(\Lambda_j\eins_{p-1}-\Lambda_{\neq j})}\Lambda_j^{n_1}\nonumber\\
 &&\times b_l^{i-1}b_m^{n_2-p}(b_m+\Lambda_j)^{p-n_1-n_2-1}\label{kernel}
\end{eqnarray}
with
\begin{align}\label{element}
e_k^{p}(\Lambda)=\sum_{1\leq j_1<\ldots<j_k\leq p} \Lambda_{j_1}\ldots\Lambda_{j_k}
\end{align}
the elementary polynomial of degree $k$ with $p$ arguments and $\Lambda_{\neq j}$  the set of $(p-1)$ eigenvalues $\Lambda$ where $\Lambda_j$ is excluded. Deriving the determinantal point process~\eqref{k-point} for the complex case is our second main result we have been aiming at.

\begin{figure}[t!]
 \centering
 \includegraphics[width=0.8\textwidth]{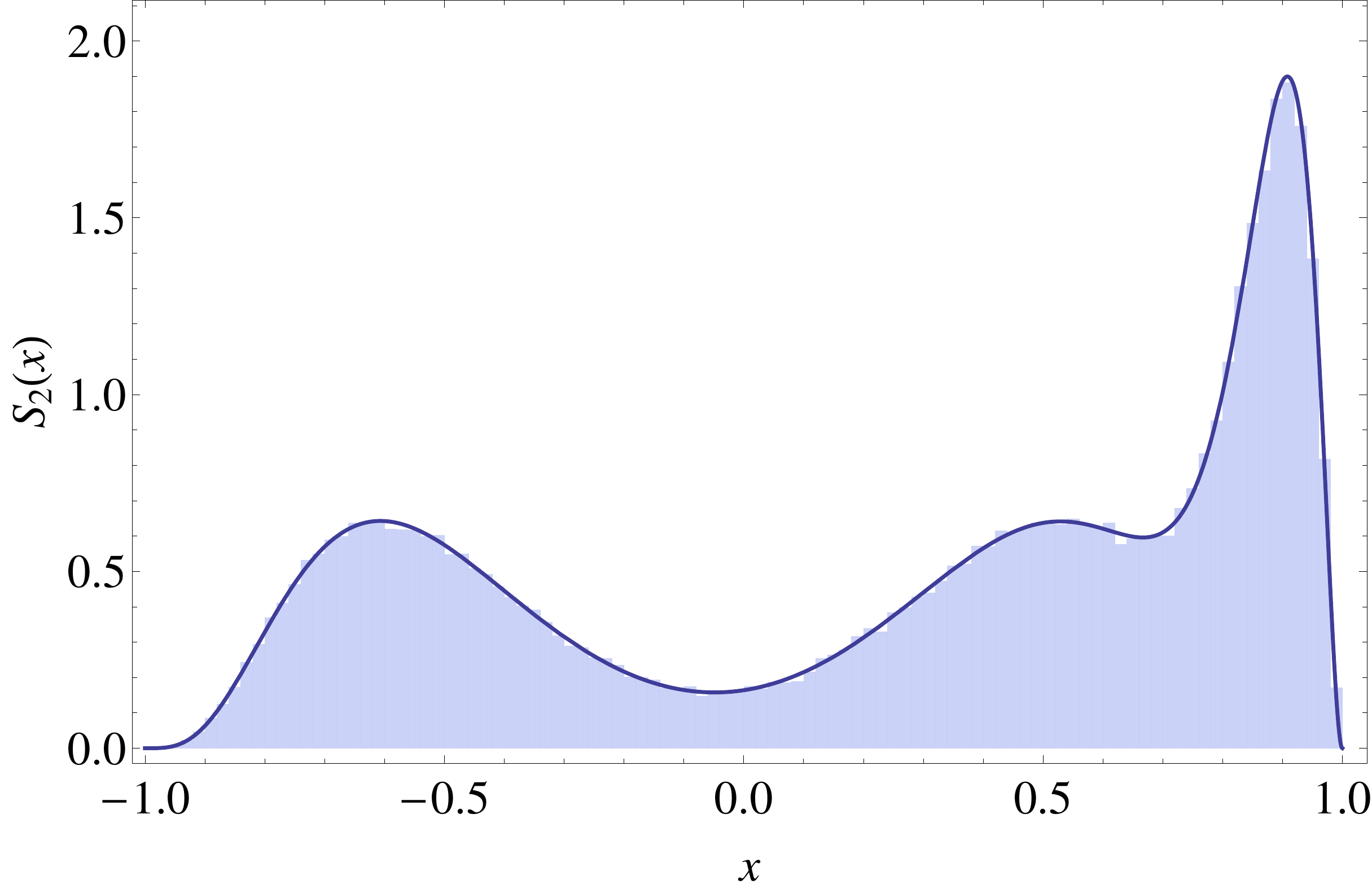}
 \caption{Comparison of the analytical expression~\eqref{lev-res} (solid curve) with Monte Carlo simulations (histogram). The random matrices have the sizes $p=3,n_1=5,n_2=7$ and the empirical eigenvalues are $\Lambda={\diag}(1/3,2,9/2)$. We have generated $50~000$ complex correlated Jacobi matrices.}
 \label{fig:JoacoiComplexDensity}
\end{figure}

The result~\eqref{k-point} immediately yields the $k$-point correlation function of the correlated Jacobi ensemble via the explicit relation
\begin{equation}\label{rel-k-point}
 \widehat{R}_{p,2}^{k}(x)=\prod_{j=1}^k\frac{2}{(1+x_j)^2}{R'}_{p,2}^k\left(\frac{1-x}{1+x}\right).
\end{equation}
The prefactor is the Jacobian from the change of coordinates $b=(1-x)/(1+x)$, cf. Eq.~\eqref{eq:k-point-CL-self}. In particular the level density~\eqref{density-def} reads
\begin{eqnarray}
 S_2(x)&=&\frac{2}{p}\sum_{i,j=1}^p (-1)^{p-i}\frac{(n_1+n_2-p)!}{(n_2-p+i-1)!(n_1-i)!} \frac{e_{p-i}^{p-1}(\Lambda_{\neq j})}{\det(\Lambda_j\eins_{p-1}-\Lambda_{\neq j})}\Lambda_j^{n_1}\nonumber\\
 &&\times \frac{1}{(1+x)^2}\left(\frac{1-x}{1+x}\right)^{n_2-p+i-1}\left(\frac{1-x}{1+x}+\Lambda_j\right)^{p-n_1-n_2-1}\label{lev-res}
\end{eqnarray}
This closed-form expression for the level density is plotted in Fig.~\ref{fig:JoacoiComplexDensity} where we compare the result with a Monte-Carlo simulation. The structure of Eq.~\eqref{lev-res} is reminiscent to the result found for the correlated complex Wishart ensemble~\cite{ATLV,SimMou}.

Despite the fact that the kernel~\eqref{kernel} is very explicit it is unsuitable for studying its asymptotics, see section~\ref{sec:saddlepoint}. This is the reason why we want to rewrite this expression in the next subsection.

\subsection{Supersymmetry and the Kernel}\label{sec:kernel}

Let us point out that the kernel~\eqref{kernel} is independent of $k$. This is quite convenient since we can express the kernel as an integral which is very similar to the one we started from, see Eqs.~\eqref{jpd} and \eqref{k-point-def-prime}. Thus the kernel is
\begin{equation}
 K(x_1,x_2;\Lambda)= K_{\rm jpd} {\det}^{n_1}\Lambda\prod_{j=1}^{p-1}\int_0^\infty db_j \frac{\Delta_p(b,x_1)}{\Delta_p(\Lambda)}\det\left[\begin{array}{c} b_i^{n_2-p}(b_i+\Lambda_j)^{p-n_1-n_2-1} \\ x_2^{n_2-p}(x_2+\Lambda_j)^{p-n_1-n_2-1}  \end{array}\right]_{\substack{1\leq i\leq p-1\\ 1\leq j \leq p}}.\label{ker1}
\end{equation}
The normalization constant is fixed via the asymptotics of the kernel in the variables $x_1$ and $x_2$.

Note that we now integrate only over $p-1$ variables instead of $p$. The missing integral can be introduced by a Dirac delta function $\delta(b_p-x_2)=\lim_{\varepsilon\to0}{\rm Im}\, 1/[\pi(b_p-x_2-\imath\varepsilon)]$. The symmetrization in all $b_1,\ldots,b_p$ reads
\begin{eqnarray}
 K(x_1,x_2;\Lambda)&=& \frac{1}{\pi p}\lim_{\varepsilon\to0}{\rm Im}\,\frac{1}{x_2-x_1}K_{\rm jpd} {\det}^{n_1}\Lambda\prod_{j=1}^{p}\int_0^\infty db_j \frac{\Delta_p(b)}{\Delta_p(\Lambda)}\frac{\det(b-x_1\eins_p)}{\det(b-(x_2+\imath\varepsilon)\eins_p)}\nonumber\\
 &&\times\det\left[\begin{array}{c} b_i^{n_2-p}(b_i+\Lambda_j)^{p-n_1-n_2-1}  \end{array}\right]_{\substack{1\leq i,j\leq p}}.\label{ker2}
\end{eqnarray}

What did we gain from rewriting the kernel? We can now identify the integral on the right hand side with a partition function~\eqref{eq:cje:lorentziangenfuncII}, i.e.
\begin{eqnarray}
 K(x_1,x_2;\Lambda)&=& \frac{1}{\pi p}\lim_{\varepsilon\to0}{\rm Im}\,\frac{1}{x_2-x_1}{Z'}_{p,2}^{1|1}(x_2+\imath\varepsilon,x_1).\label{ker3}
\end{eqnarray}
This allows us to apply the result~\eqref{eq:cje:chifinal}  of the projection formula, 
\begin{eqnarray}
 K(x_1,x_2;\Lambda)&=&\lim_{\varepsilon\to0}\frac{1}{2\pi\imath p}\sum_{L=\pm1}\frac{L}{x_2-x_1}\int d\mu(U) \sdet^{-1}(\eins_p\otimes U-\Lambda^{-1}\otimes{\diag}(x_2+\imath L\varepsilon,x_1))\nonumber\\
&&\times
\sdet^{p-n_1-n_2}(L U+\eins_{1|1})\sdet^{n_2}U.
\label{projection-res}
\end{eqnarray}
The explicit parametrization of $U$ is
\begin{equation}
U=\left[\begin{array}{cc} L e^\vartheta & \eta^* \\ \eta & e^{\imath\varphi} \end{array}\right]
\end{equation}
with $\eta$ and $\eta^*$ two independent Grassmann variables and $\vartheta\in\mathbb{R}$ and $\varphi\in[0,2\pi]$. The Haar measure is $d\mu(U)=Le^{\vartheta+\imath\varphi}d\vartheta d\varphi d\eta^* d\eta/(2\pi)$ where we choose the  convention $\int \eta d\eta=1$. The expression~\eqref{projection-res} is exactly what we are aiming at. The integral is over a small fixed number of variables and the ``large'' dimensions $n_1$, $n_2$ and $p$ appear as external parameters. Thus the result~\eqref{projection-res} invites for a saddle point approximation. This is the main idea behind the supersymmetry method~\cite{Efebook,Guhrbook}.

The diagonal elements of the determinant~\eqref{k-point} are the level density at positions $x_1,\ldots,x_k$ whereas the kernel becomes
\begin{eqnarray}
S'_{2}(b)&=&\lim_{\varepsilon\to0}\sum_{L=\pm1}\frac{L}{4\pi\imath p}\int d\mu(U) \str(\eins_p\otimes U-(b+\imath L\varepsilon)\Lambda^{-1}\otimes\eins_{1|1})^{-1}\left[\begin{array}{cc} \eins_p & 0 \\ 0 & -\eins_p \end{array}\right]\nonumber\\
&&\times\sdet^{-1}(\eins_p\otimes U-(b+\imath L\varepsilon)\Lambda^{-1}\otimes\eins_{1|1})
\sdet^{p-n_1-n_2}(L U+\eins_{1|1})\sdet^{n_2}U\nonumber\\
\label{projection-lev}
\end{eqnarray}
due to L'H\^{o}spital's rule.

\section{Eigenvalue Density in the Real Ensemble}\label{sec:cje:realcase}

In the case of the real correlated Jacobi ensemble the computation of an arbitrary $k$-point correlation function is highly non-trivial since the corresponding group integral~\eqref{group-int} is unknown. Hence we concentrate on the calculation of the level density~\eqref{density-def}. Already computing the level density with the help of the projection formula is quite involved. The dimension of the supermatrix model is $(2|2)\times(2|2)$ and thus twice as large as for the kernel of the complex matrices.

We apply
the generalized Hubbard-Stratonovich transformation~\cite{IngSie,KGG09,KSG} for $k=1$, third equality of Eq.~(\ref{eq:cje:chifinal}), which explicitly reads in this case
\begin{eqnarray}
{Z'}_{p,1}^{1}(\kappa) &=& \frac{1}{(4\pi)^2}\int d[\sigma]I_{n_2}(\sigma)\sdet^{-\mu/2}\left(\sigma + \eins_{2|2}\right)\sdet^{-1/ 2}\left(\Lambda^{-1}\otimes\kappa -
\eins_p\otimes\sigma\right).
\label{eq:cje:Rdensity:HS}
\end{eqnarray}
The exponent is $\mu=n_1+n_2-p$ and $I_{n_2}(\sigma)$ is the supersymmetric Ingham-Siegel integral~\eqref{eq:gHS:InghamSiegel} for $k_1=k_2=1$. Furthermore we have $\kappa={{\diag}}( b\pm\imath  \varepsilon,b\pm\imath  \varepsilon,b_1,b_1)$
and
\begin{align}
 \sigma = \begin{pmatrix} \sigma_1 & \zeta\\-\zeta^\dagger & i\sigma_2\eins_2\end{pmatrix}~, \quad \zeta = \begin{pmatrix} \alpha & \alpha^* \\ \beta & \beta^*\end{pmatrix}\quad \sigma_1 =
 \begin{pmatrix} a & b\\ b& c\end{pmatrix}.\label{eq:cje:Rdensity:parametrization}
\end{align}
The normalization has been calculated by choosing $b_1=b\pm\imath  \varepsilon$  where Cauchy-like integral theorems~\cite{Wegner,ParSou,Efetov,Con,ConGro,KieburgKohlerGuhr} apply. In the parametrization of $\sigma$,  the variables $a,b,c,\sigma_2$ are real commuting whereas $\alpha,\beta$ are complex Grassmann variables. The measure is
\begin{align}
d[\sigma]=dad{b}d{c}d{\sigma_2}d{\alpha}d{\alpha^*}d{\beta}d{\beta^*}.
\end{align}

Let us first consider the diagonal blocks. Since the boson-boson as well as the fermion-fermion block of $\kappa$ are proportional to $\eins_2$, we can diagonalize the boson-boson block of $\sigma$ without any problems. Accordingly, we write $\sigma_1=OrO^T$, where $r={{\diag}}(r_1,r_2)$ and  $O\in\text{O}(2)$. This change of coordinates yields a decomposition of the differential
\begin{align}
 d[\sigma_1] = \left|r_1-r_2\right|dr_1 dr_2 d{\mu(O)}.
\end{align}
Because of the structure of the integrand~(\ref{eq:cje:Rdensity:HS}), the integral over ${\rm O}(2)$ factorizes and yields a factor of $\pi/2$.

In the next step, we expand the superdeterminants of $\sigma+\eins_{2|2}$ in the Grassmann variables $\alpha,\alpha,\beta$ and $\beta^*$,
\begin{eqnarray}
&&\sdet^{-\mu/2}\left(\sigma + \eins_{2|2}\right)= \frac{(1+\imath \sigma_2)^{\mu}}{(1+r_1)^{\mu/2}(1+r_2)^{\mu/2}}\label{eq:cje:Rdensity:npdetexp}\\
&&\times\left(1+ \frac{\mu}{(1+\imath \sigma_2)
 (1+r_1)}\alpha^*\alpha + \frac{\mu}{(1+\imath \sigma_2)(1+r_2)}\beta^*\beta+\frac{\mu(\mu-1)}{(1+\imath \sigma_2)^2(1+r_1)(1+r_2)}\alpha^*\alpha \beta^*\beta\right)\nonumber
\end{eqnarray}
and
\begin{eqnarray}
  &&\hspace*{-1.5cm}\sdet^{-1/ 2}\left(\Lambda^{-1}\otimes\kappa -
\eins_p\otimes\sigma\right)=\frac{\det(\kappa_2\Lambda^{-1}-\imath \sigma_2\eins_p)}{\sqrt{\det(\kappa_1\Lambda^{-1}-r_1\eins_p)\det(\kappa_1\Lambda^{-1}-r_2\eins_p)}}\nonumber \\
  &&\times\biggl(1+\sum_{j=1}^p\frac{1}{(\kappa_2\Lambda_j^{-1}-\imath \sigma_2)(\kappa_1\Lambda_j^{-1}-r_1)}\alpha^*\alpha+\sum_{j=1}^p\frac{1}{(\kappa_2\Lambda_j^{-1}-\imath \sigma_2)(\kappa_1\Lambda_j^{-1}-r_2)}\beta^*\beta\nonumber\\
  &&+\sum_{1\leq i\neq j\leq p}\frac{1}{(\kappa_2\Lambda_i^{-1}-\imath \sigma_2)(\kappa_2\Lambda_j^{-1}-\imath \sigma_2)(\kappa_1\Lambda_i^{-1}-r_1)(\kappa_1\Lambda_j^{-1}-r_2)}\alpha^*\alpha\beta^*\beta\biggl).\label{sdet-Lambda}
\end{eqnarray}
The expansion in the Grassmann variables of the supersymmetric Ingham-Siegel integral for $\beta=1$ was done in Ref.~\cite{Rech12} and is
\begin{eqnarray}
 I_{n_2}(r)&=&\frac{4\pi}{(n_2-2)!}\Theta(r_1)\Theta(r_2)(r_1r_2)^{(n_2-1)/2}\left(\imath\frac{\partial}{\partial\sigma_2}\right)^{n_2-2}\nonumber\\
&&\times\left(1 - \imath\left(\frac{\alpha^*\alpha}{r_1}+\frac{\beta^*\beta}{r_2}\right)\frac{\partial}{\partial \sigma_2}
-\frac{\alpha^*\alpha\beta^*\beta}{r_1r_2}\frac{\partial^2}{\partial \sigma_2^2}\right)\delta(\sigma_2),
\label{eq:cje:Rdensity:InghamSiegelExp}
\end{eqnarray}
where $\Theta$ is the Heaviside step function.

When performing the integrals over the Grassmann variables we only keep the leading order terms in the four Grassmann variables, in particular those terms proportional to $\alpha^*\alpha\beta^*\beta$. Moreover we evaluate the Dirac delta function and apply the derivative in $\kappa_2$. In the end we set $\kappa_2=\kappa_1=\kappa$.
Then the partition function is
\begin{eqnarray}
\frac{\partial {Z'}_{p,1}^{1}}{\partial \kappa_2}(\kappa,\kappa)&=&\frac{1}{8}\left[\mu(\mu-1)C_{n_2-2,\mu-2}^p(\kappa;\Lambda)N_{n_2-1,\mu,+1,+1}^0(\kappa;\Lambda)-2\mu(n_2-1)C_{n_2-1,\mu-1}^p(\kappa;\Lambda)\right.\nonumber\\
&&\hspace*{-3cm}\left.\times N_{n_2-1,\mu,+1,-1}^0(\kappa;\Lambda)+n_2(n_2-1)C_{n_2,\mu}^p(\kappa;\Lambda)N_{n_2-1,\mu,-1,-1}^0(\kappa;\Lambda)\right]\nonumber\\
&&\hspace*{-3cm}+\frac{1}{4}\sum_{j=1}^p\left[\mu C_{n_2-2,\mu-1}^{p-1}(\kappa;\Lambda_{\neq j})N_{n_2-1,\mu,+1}^1(\kappa;\Lambda;\Lambda_j)-(n_2-1)C_{n_2-1,\mu}^{p-1}(\kappa;\Lambda_{\neq j}) N_{n_2-1,\mu,-1}^1(\kappa;\Lambda;\Lambda_j)\right]\nonumber\\
&&\hspace*{-3cm}+\frac{1}{4}\sum_{1\leq i<j\leq p}C_{n_2-2,\mu}^{p-2}(\kappa;\Lambda_{\neq i,j})N_{n_2-1,\mu}^2(\kappa;\Lambda;\Lambda_i,\Lambda_j)\label{eq:cje:Rdensity:derZ}
\end{eqnarray}
with $\Lambda_{\neq j}$ the diagonal $p-1$ dimensional sub-matrix of $\Lambda$ removing $\Lambda_j$ and $\Lambda_{\neq i,j}$ the diagonal $p-2$ dimensional sub-matrix removing $\Lambda_i$ and $\Lambda_j$. The functions $C_{a,b}^c(\kappa;E_1,\ldots,E_c)$ are the integrals over $\sigma_2$. The derivatives of the Dirac delta functions can be written as contour integrals and as a finite sum
\begin{eqnarray}
 C_{a,b}^c(\kappa;E_1,\ldots,E_c)&=&\frac{\partial}{\partial\kappa}\int_0^{2\pi} \frac{d\varphi}{2\pi} e^{-\imath a \varphi}(1+e^{\imath\varphi})^b\prod_{j=1}^c(\kappa E_j^{-1}-e^{\imath\varphi})\nonumber\\
  &=&\sum_{j=1}^{c}(-1)^{c-j}\frac{b!\ j}{(c-a-j)!(a+b-c+j)!} \kappa^{j-1} e^c_{j}(E^{-1}).\label{int0}
\end{eqnarray}
The resulting sum involves the elementary polynomials~\eqref{element} which already appeared for the correlated complex ensembles. The function
\begin{equation}
 N_{a,b,d_1,d_2}^0(\kappa;\Lambda)=\int_{\mathbb{R}_+^2}\frac{dr_1dr_2|r_1-r_2|}{\prod_{l=1}^p\sqrt{\kappa\Lambda_l^{-1}-r_1}\sqrt{\kappa\Lambda_l^{-1}-r_2}}\frac{r_1^{(a-1+d_1)/2}r_2^{(a-1+d_2)/2}}{(1+r_1)^{(b+1+d_1)/2}(1+r_2)^{(b+1+d_2)/2}}\label{int1}
\end{equation}
appears for any term in the expansion of the Grassmann variables as long as we take the first term in the superdeterminant~\eqref{sdet-Lambda}. The second and third term in the expansion~\eqref{sdet-Lambda} yields the two-fold integrals
\begin{equation}
 N_{a,b,d}^1(\kappa;\Lambda;\Lambda_j)=\int_{\mathbb{R}_+^2}\frac{dr_1dr_2|r_1-r_2|}{\prod_{l=1}^p\sqrt{\kappa\Lambda_l^{-1}-r_1}\sqrt{\kappa\Lambda_l^{-1}-r_2}}\frac{1}{\kappa \Lambda_j^{-1}-r_1}\frac{r_1^{a/2}r_2^{(a-1+d)/2}}{(1+r_1)^{b/2}(1+r_2)^{(b+1+d)/2}}\label{int2}
\end{equation}
and 
\begin{eqnarray}
 N_{a,b}^2(\kappa;\Lambda;\Lambda_i,\Lambda_j)&=&\int_{\mathbb{R}_+^2}\frac{dr_1dr_2|r_1-r_2|}{\prod_{l=1}^p\sqrt{\kappa\Lambda_l^{-1}-r_1}\sqrt{\kappa\Lambda_l^{-1}-r_2}}\frac{1}{(\kappa \Lambda_i^{-1}-r_1)(\kappa\Lambda_j^{-1}-r_2)}\nonumber\\
 &&\times\frac{r_1^{a/2}r_2^{a/2}}{(1+r_1)^{b/2}(1+r_2)^{b/2}},\label{int3}
\end{eqnarray}
respectively. These two integrals are principal value integrals at the non-integrable singularities $r_1=\kappa \Lambda_j^{-1}$ and $r_2=\kappa \Lambda_i^{-1}$. We discuss the numerical evaluation of these two integrals in appendix~\ref{app:expl}.

The structure of Eq.~(\ref{eq:cje:Rdensity:derZ}) is very similar to that obtained for the
ordinary and the doubly correlated Wishart model computed in Refs.~\cite{Rech12,WaltnerWirtzGuhr2014}. The full expression can be separated into three parts; a part with only square root singularities in $r_1$ and $r_2$, a part including an additional $3/2$-singularity in $r_1$, and a part with two $3/2$-singularities, one for $r_1$ and one for $r_2$. The latter two terms in Eq.~\eqref{eq:cje:Rdensity:derZ} correspond to the single and double sum. They have to be regularized by Cauchy principal value integrals which is done in appendix~\ref{app:expl}.

The $\varepsilon\rightarrow0$ limit of the imaginary part is taken from the integrals~\eqref{int1}, \eqref{int2}, and \eqref{int3}, only. Thereby we adapt the analysis from Refs.~\cite{Rech12,WaltnerWirtzGuhr2014} and concentrate on the product of the square roots. We choose the branch cut of each square root along the negative real line. Then for any $y\in \mathbb{R}$  we have
\begin{align}
 \lim_{\varepsilon\to 0} \frac{1}{\sqrt{y-\imath L\varepsilon}}=
  \frac{\Theta(y)+L\imath\Theta(-y)}{\sqrt{|y|}}.\label{eq:cje:Rdensity:exampleIII}
\end{align}
Recall that $L$ is the sign of the imaginary increment. Thus the imaginary part of Eq.~(\ref{eq:cje:Rdensity:derZ}) is only non-vanishing in the limit $\varepsilon\rightarrow0$ if $\det(b\Lambda^{-1}-r_1\eins_p)\det(b\Lambda^{-1}-r_2\eins_p)$ has a negative real part. We recall that $\kappa=b+\imath\varepsilon$. Let us assume that $0<\Lambda_1<\Lambda_2<\cdots<\Lambda_p<\infty$. Then we can divide the integration domain $[0,\infty)$ into disjoint subsets
\begin{align}
 [0,\infty) = \bigcup_{i=0}^p \text{V}_i
\end{align}
where
\begin{align}\label{domain}
\text{V}_0=\left[0,b\Lambda_p^{-1}\right),\ \text{V}_p=\left(b\Lambda_1^{-1},\infty\right),\ {\rm and}\ \text{V}_i=\left(b\Lambda_{p-i+1}^{-1},b\Lambda_{p-i}^{-1}\right)\ \text{for}\ i=1,\dots,p-1.
\end{align}
This decomposition implies $\det(b\Lambda^{-1}-r\eins_p)=(-1)^{l}|\det(b\Lambda^{-1}-r\eins_p)|$ for $r\in \text{V}_l$.

We plug the decomposition~\eqref{domain} into the integrals~\eqref{int1}, \eqref{int2}, and \eqref{int3} and keep only those terms which  yield an imaginary part. Then the double integral becomes a sum of decoupled one-fold integrals
\begin{equation}\label{int-sum}
 \int_0^\infty dr_1\int_0^\infty dr_2|r_1-r_2|\rightarrow\sum_{\substack{0\leq l_1,l_2\leq p \\ l_1+l_2\in2\mathbb{N}_0+1}}{\rm sign}(l_1-l_2)\int_{V_{l_1}}dr_1\int_{V_{l_2}}dr_2(r_1-r_2).
\end{equation}
Moreover, we get an additional sign in the sum from the square roots. Assuming $r_1\in\text{V}_{l_1}$ and $r_2\in\text{V}_{l_2}$, we have
\begin{align}
\begin{split}
&\text{Im}\lim_{\varepsilon\rightarrow0}\prod_{k=1}^p\frac{1}{\sqrt{\left((b+\imath\varepsilon)\Lambda_k^{-1}-r_1\right)}\sqrt{\left((b+\imath\varepsilon)\Lambda_k^{-1}-r_2\right)}}
\\&=\left\{\begin{array}{cl}\displaystyle\prod_{k=1}^p\frac{(-1)^{(l_1+l_2+1)/2}}{\sqrt{\left|b\Lambda_k^{-1}-r_1\right|\left|b\Lambda_k^{-1}-r_2\right|}}, & l_1+l_2\in 2\mathbb{N}+1,\\0, & \text{otherwise}.\end{array}\right.\end{split}\label{eq:cje:Rdensity:limitvarepsilon}
\end{align}
We also obtain for each additional term $1/(\kappa\Lambda_j^{-1}-r)$ in the integrals~\eqref{int2} and \eqref{int3} a sign ${\rm sign}(p-l-j+1/2)$ if $r\in\text{V}_{l}$. The term $1/2$ in the ${\rm sign}$ function guarantees that the sign is positive if $j=p-l$.

Summarizing everything the level density~\eqref{density-CL-def} of the correlated Cauchy-Lorentz ensemble is
\begin{align}
\begin{split}
 &S'_1(b)=\sum_{\substack{0\leq l_1,l_2\leq p\\l_1+l_2\in2\mathbb{N}_0+1}}\int\limits_{V_{l_1}\times V_{l_2}}d{r_1}d{r_2}\frac{(-1)^{(l_1+l_2-1)/2}}{{\displaystyle\prod_{k=1}^p}\sqrt{\left|b\Lambda_k^{-1}-r_1\right|\left|b\Lambda_k^{-1}-r_2\right|}}
 \\&\times\left(f_{1}(r_1,r_2;b,\Lambda)+\sum_{l=1}^p\frac{f_{2,l}(r_1,r_2;b,\Lambda)}{\left(b\Lambda_l^{-1}-r_1\right)}+\sum_{\substack{i,l=1\\i\neq l}}^p\frac{f_{3,l,i}(r_1,r_2;b,\Lambda)}{\left(b\Lambda_l^{-1}-r_1\right)\left(b\Lambda_i^{-1}-r_2\right)}\right),
\end{split} \label{eq:cje:Rdensity:finaldensity}
\end{align}
where $f_{1}(r_1,r_2;b,\Lambda)$, $f_{2,l}(r_1,r_2;b,\Lambda)$, and $f_{3,l,i}(r_1,r_2;b,\Lambda)$ are some functions that can be read off from Eq.\ (\ref{eq:cje:Rdensity:derZ}). As already said the integrals with the $3/2$ singularities have to be regularized which is done in appendix~\ref{app:expl}.
Finally we obtain a finite sum where each summand is a product of the three integrals $C_{a,b}^c(\kappa;E_1,\ldots,E_c)$ (see Eq.~\eqref{int0}), $g_{a,c,l}^0(b;\Lambda)$ (see Eq.~\eqref{int-g1}), and $g_{a,c,l}^1(b;\Lambda;\Lambda_i)$ (see Eq.~\eqref{int-g2}),
\begin{eqnarray}
S'_1(b)&=&\frac{1}{8\pi}\sum_{\substack{0\leq l_1,l_2\leq p \\ l_1+l_2\in2\mathbb{N}_0+1}}(-1)^{(l_1+l_2+1)/2}{\rm sign}(l_1-l_2)\label{CL-dens-result}\\
&&\hspace*{-1cm}\times\biggl(\mu\det\left[\begin{array}{ccc} 0 & g_{n_2+1,\mu+2,l_1}^0(b;\Lambda) & g_{n_2-1,\mu+2,l_1}^0(b;\Lambda) \\ (n_2-1) C_{n_2-1,\mu-1}^{p}(b;\Lambda) & g_{n_2+1,\mu+2,l_2}^0(b;\Lambda) & g_{n_2-1,\mu+2,l_2}^0(b;\Lambda) \\ (\mu-1) C_{n_2-2,\mu-2}^{p}(b;\Lambda) & g_{n_2-1,\mu,l_2}^0(b;\Lambda) & g_{n_2-3,\mu,l_2}^0(b;\Lambda) \end{array}\right]\nonumber\\
&&\hspace*{-1cm}+(n_2-1)\det\left[\begin{array}{ccc} 0 & g_{n_2-1,\mu,l_1}^0(b;\Lambda) & g_{n_2-3,\mu,l_1}^0(b;\Lambda) \\ \mu C_{n_2-1,\mu-1}^{p}(b;\Lambda) & g_{n_2-1,\mu,l_2}^0(b;\Lambda) & g_{n_2-3,\mu,l_2}^0(b;\Lambda) \\ n_2 C_{n_2,\mu}^{p}(b;\Lambda) & g_{n_2+1,\mu+2,l_2}^0(b;\Lambda) & g_{n_2-1,\mu+2,l_2}^0(b;\Lambda) \end{array}\right]\nonumber\\
&&\hspace*{-1cm}+2\sum_{j=1}^p\det\left[\begin{array}{ccc} 0 & g_{n_2+1,\mu,l_1}^1(b;\Lambda;\Lambda_j) & g_{n_2-1,\mu,l_1}^1(b;\Lambda;\Lambda_j) \\ (n_2-1)C_{n_2-1,\mu}^{p-1}(b;\Lambda_{\neq j}) & g_{n_2+1,\mu+2,l_2}^0(b;\Lambda) & g_{n_2-1,\mu+2,l_2}^0(b;\Lambda) \\ \mu C_{n_2-2,\mu-1}^{p-1}(b;\Lambda_{\neq j}) & g_{n_2-1,\mu,l_2}^0(b;\Lambda) & g_{n_2-3,\mu,l_2}^0(b;\Lambda) \end{array}\right]\nonumber\\
&&\hspace*{-1cm}+2\sum_{1\leq i<j\leq p}C_{n_2-2,\mu}^{p-2}(b;\Lambda_{\neq i,j})\det\left[\begin{array}{cc} g_{n_2+1,\mu,l_1}^1(b;\Lambda;\Lambda_i) & g_{n_2-1,\mu,l_1}^1(b;\Lambda;\Lambda_i) \\ g_{n_2+1,\mu,l_2}^1(b;\Lambda;\Lambda_j) & g_{n_2-1,\mu,l_2}^1(b;\Lambda;\Lambda_j) \end{array}\right]\biggl).\nonumber
\end{eqnarray}
The integrals $g_{a,c,l}^0(b;\Lambda)$ and $g_{a,c,l}^1(b;\Lambda;\Lambda_i)$ are one-fold integrals over the compact interval $V_l$ which is numerically more advantageous than the original two-fold integral.

\begin{figure}[t!]
\centering
\includegraphics[width=0.8\textwidth]{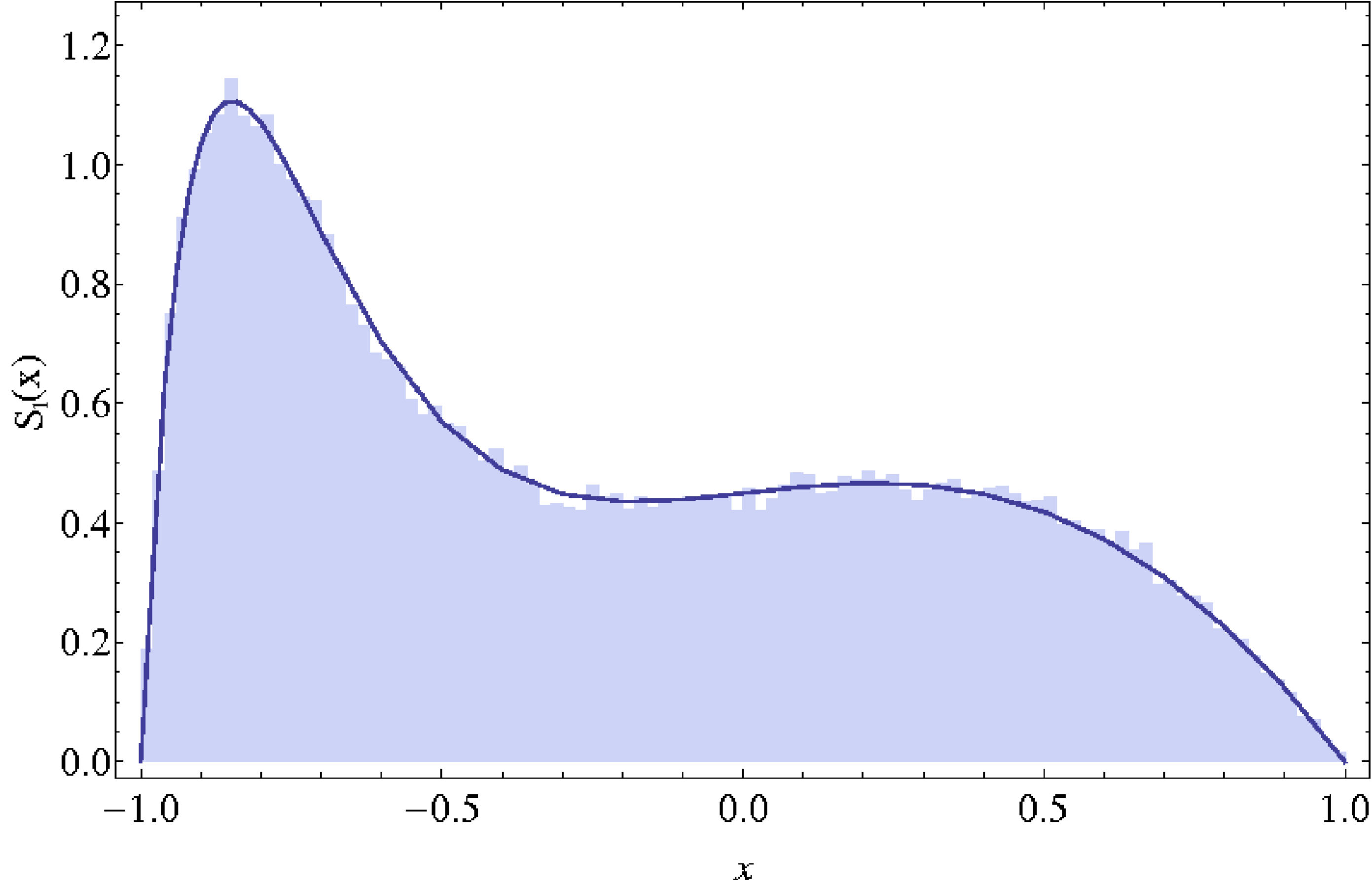}
\caption{Comparison between the result~(\ref{eq:cje:Rdensity:finaldensity}) (solid curve) and Monte Carlo simulations (histogram) for $n_1=5$, $n_2=5$ and $p=2$ with $\Lambda_j=1,\,4$. The sample size consists of $50~000$ real correlated Jacobi matrices.}
\label{fig:realdensity}
\end{figure}

 The level density~\eqref{density-def} of the correlated Jacobi ensemble readily follows from~\eqref{CL-dens-result} via the relation
\begin{align}
S_1(x)=\frac{2}{(1+x)^2}S'_1\left(\frac{1-x}{1+x}\right).\label{eq:cje:Rdensity:finaldensity}
\end{align}
This third main result is compared to Monte Carlo simulations in Fig.~\ref{fig:realdensity}. The excellent agreement validates our calculation.

We emphasize that the integrals drastically simplify when the spectrum of $C_{\rm eff}$ and thus $\Lambda$ is doubly degenerate. Then the limit $\varepsilon\to0$ in combination with taking the imaginary part yields Dirac delta functions such that one of the integrals over $r_1$ and $r_2$ can be exactly performed. The remaining integral can be easily performed, too, by employing the idea of~\cite{KieGuh} where the term $1/\det(b \Lambda-r\eins_p)$ can be written as a quotient of a Vandermonde determinant and a Cauchy-Vandermonde determinant. The degeneracy of the spectrum is not that academic. Such a degeneracy was proposed for the correlated real Wishart ensemble in \cite{WirtzKieburgGuhrPRL} to calculate explicit analytical results. It was shown that in the limit of large matrix dimensions such an artificially introduced degeneracy has no effect on the spectral properties like the level density and the $k$-point correlation function. The open question is if this statement carries over to the Cauchy-Lorentz 
and the Jacobi ensemble. For the level density we answer this question in the next section.

\section{Asymptotics of the Level Density}\label{sec:saddlepoint}

To obtain the asymptotic behavior of the level density for both, the real and the complex correlated Jacobi models, we perform a saddle point approximation of the supersymmetric expression~\eqref{eq:cje:chifinal} for $k=1$. Particularly we first consider the level density of the correlated Cauchy-Lorentz ensemble
\begin{eqnarray}
S'_\beta(b) &=&\frac{1}{ \pi p }\lim_{\varepsilon\rightarrow0}\text{Im}\left.\frac{\partial}{\partial J}\right|_{J=0}{Z'}_{p,\beta}^{1}(b+\imath \varepsilon,b-J+\imath \varepsilon)\nonumber\\
&=& \frac{K}{\pi p\gamma}\lim_{\varepsilon\rightarrow0}\text{Im} \int d\mu(U) \exp\left[-\frac{p}{\gamma}\mathcal{L}(U)\right]\sum_{j=1}^p\str(b+\imath\varepsilon-\Lambda_j U)^{-1}\left[\begin{array}{cc} \eins_\gamma & 0 \\ 0 & 0 \end{array}\right].
\label{densitySaddlePoint}
\end{eqnarray}
The normalization constant $K$ is
\begin{equation}\label{norm}
K^{-1}= \int d\mu(U) \exp\left[-\frac{p}{\gamma}\mathcal{L}(U)\right],
\end{equation}
where we used an apparent $b$-dependent version which is more convenient in the saddle point approximation. The $b$-dependence is indeed only apparent because of the Cauchy-like integration theorems in superspaces~\cite{Wegner,ParSou,Efetov,Con,ConGro,KieburgKohlerGuhr}.

 The ``Lagrangian'' in the exponential function is 
\begin{align}
 \mathcal{L}(U) =\str\left(\frac{n_1+n_2-p}{p}\ln\left(U + \eins_{ \gamma| \gamma}\right)- \frac{n_2}{p}\ln U+\frac{1}{p}\sum_{k=1}^p\ln\left( (b+\imath\varepsilon)\Lambda_{k}^{-1}\eins_{ \gamma| \gamma} -U\right)\right).\label{lagrangian}
\end{align}
This Lagrangian has to be minimized if $n_1\propto n_2\propto p\gg1$, meaning that its first derivative has to vanish,
\begin{equation}
\mathcal{L}'(U_0)=\frac{n_1+n_2-p}{p}(U_0+\eins_{ \gamma| \gamma})^{-1}- \frac{n_2}{p}U_0^{-1}-\frac{1}{p}\sum_{k=1}^p\left( b\Lambda_{k}^{-1}\eins_{ \gamma| \gamma} -U_0\right)^{-1}=0.
\end{equation}
Since this equation is invariant under the supergroup ${\rm UOSp}(2|2)$ for $\beta=1$ and ${\rm U}(1|1)$ for $\beta=2$ we may diagonalize $U_0$ such that each of its eigenvalues satisfies
\begin{equation}\label{saddle}
\mathcal{L}'(q_0)=\frac{n_1+n_2-p}{p}\frac{1}{q_0+1}- \frac{n_2}{p}\frac{1}{q_0}-\frac{1}{p}\sum_{j=1}^p\frac{\Lambda_j}{b -q_0\Lambda_j}=0.
\end{equation}
We perform the same asymptotic analysis as in \cite{AKI,WirtzKieburgGuhrPRL} and count $p+2$ poles at $q_0=-1,0,b\Lambda_p^{-1}, \ldots,b\Lambda_1^{-1}$. The asymptotic behavior of $\mathcal{L}'(q_0)$ at $q_0=b\Lambda_p^{-1}, \ldots,b\Lambda_1^{-1}$ implies that $p-1$ of the $p+1$ solutions of the saddle point equation~\eqref{saddle} are real.  When taking the imaginary part in Eq.~\eqref{densitySaddlePoint} we recognize that the real solutions do not contribute. Thus we are looking for the complex conjugate pair which solves Eq.~\eqref{saddle}.

To find a closed form of the saddle point solution we underline that all poles and all zero points, apart from the complex conjugate pair, of $\mathcal{L}'(q_0)$ lie on the real line. Thus an integral along an appropriate contour $C_R$ of its logarithmic derivative yields the solution via Cauchy's integration theorem, i.e.
\begin{eqnarray}
 q_0=\lim_{R\to\infty}\int_{C_R}z\frac{\partial\ln\mathcal{L}'}{\partial z}(z)\frac{dz}{2\pi\imath}.\label{solution}
\end{eqnarray}
The contour is chosen as follows
\begin{equation}
 C_R=\left\{\left.R e^{\imath\varphi}+\frac{\imath}{R}\right|\varphi\in[0,\pi]\right\}\cup\left\{\left.r+\frac{\imath}{R}\right|r\in[-R,R]\right\}
\end{equation} 
and we integrate clockwise. In this way we obtain either the unique complex solution in the complex upper half-plane or zero. The latter implies $S'_\beta(b)=0$, too. 

The fermion-fermion block encircles the pole at the origin but the other poles of the Lagrangian are purely zero points of the integrand. Hence this contour can be deformed to go through both complex solutions $q_0$ and $q_0^*$. However the boson-boson block is hindered to reach both solutions since the poles of $\mathcal{L}$ are apart from the one at origin poles of the bosonic integrand. The sign in front of the imaginary increment $\imath\varepsilon$ dictates which pole can be reached and which not.   By a contour deformation we can include the saddle point in the upper or lower complex half-plane. In Ref.~\cite{Verbaarschotetal} it was shown that the leading order in $p$ is given by those saddle points where the boson-boson block and the fermion-fermion block share the same saddle point. Thus, in the vicinity of the saddle point we have
\begin{align}
 U = \left(\text{Re}\,q_0 + L\imath  \text{Im}\,q_0\right) \left(\eins_{\tilde \gamma| \tilde \gamma} + \frac{\delta U}{\sqrt{p}}\right)
\end{align}
with $\delta U$ a Hermitian $(\gamma|\gamma)\times(\gamma|\gamma)$ supermatrix satisfying certain symmetries in the case $\beta=1$, see \cite{KGG09,KSG}. We plug this expansion into Eq.~\eqref{densitySaddlePoint} and integrate over $\delta U$. We arrive at
\begin{align}
\begin{split}
 S'_\beta(b) \approx \frac{1}{\pi}\text{Im}\frac{1}{p}\sum_{j=1}^p\frac{1}{b+\imath\varepsilon - \Lambda_j q_0(b)}=\frac{1}{\pi p}\sum_{j=1}^p\frac{\Lambda_j\text{Im}\,q_0(b)}{(b-\Lambda_j\text{Re}\,  q_0(b))^2+(\Lambda_j\text{Im}\,q_0(b))^2},
\end{split}
\label{densitySaddlePoint1}
\end{align}
where we already fixed the normalization by the integration $\int S'_\beta(b) db=1$. As in Ref.~\cite{WirtzKieburgGuhrPRL}, we can simplify the expression~\eqref{densitySaddlePoint1} with the help of the saddle point equation~\eqref{saddle}. Then we find our fourth main result for the correlated Cauchy-Lorentz ensemble
\begin{align}
\begin{split}
 S'_\beta(b) \approx \frac{n_1+n_2-p}{\pi p}\frac{1}{b}\text{Im}\frac{q_0(b)}{q_0(b)+1}=\frac{n_1+n_2-p}{\pi p}\frac{1}{b}\frac{\text{Im}\,q_0(b)}{(\text{Re}\,q_0(b)+1)^2+(\text{Im}\,q_0(b))^2}
\end{split}
\label{densitySaddlePoint2}
\end{align}
and for the correlated Jacobi ensemble
\begin{align}
\begin{split}
 S_\beta(x) \approx \frac{2(n_1+n_2-p)}{\pi p}\frac{1}{1-x^2}\frac{\text{Im}\,q_0([1-x]/[1+x])}{(\text{Re}\,q_0([1-x]/[1+x])+1)^2+(\text{Im}\,q_0([1-x]/[1+x]))^2}
\end{split}
\label{densitySaddlePoint3}
\end{align}
This solution is slightly more involved compared to the correlated Wishart ensemble~\cite{WirtzKieburgGuhrPRL}. Nevertheless it is a closed form in combination with the saddlepoint solution~\eqref{solution}.

\begin{figure}[t!]
 \centering
 \includegraphics[width=0.7\textwidth]{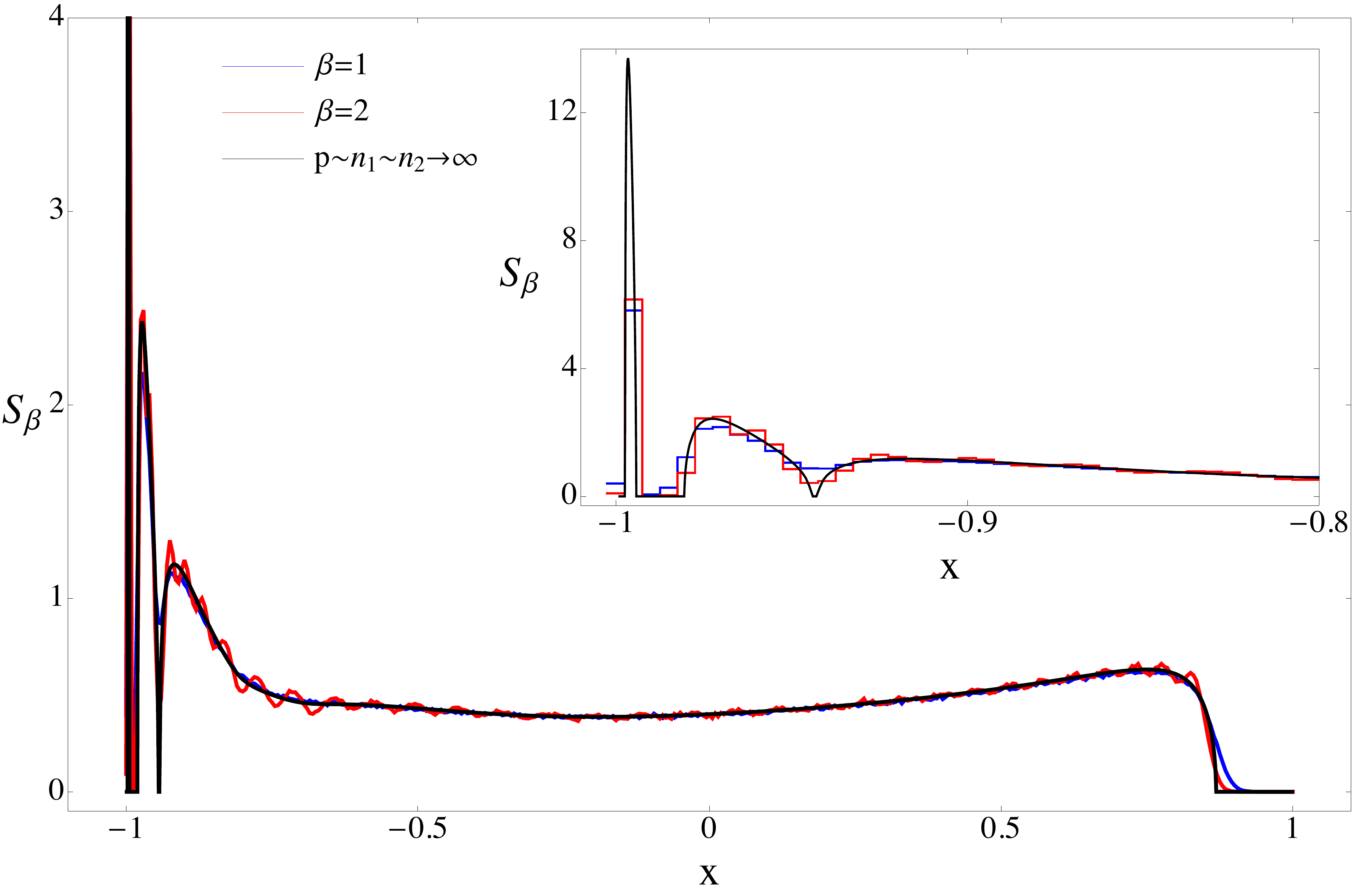}
 \caption{Comparison of the analytical expression~\eqref{densitySaddlePoint3} for the limiting eigenvalue density (solid black curve) with numerical simulations with correlated real Jacobi matrices (blue histogram) and correlated complex Jacobi matrices (red histogram). Each ensemble consists of $10^5$ matrices with the parameters  $p=32,n_1=71,n_2=68$. The empirical eigenvalues were independently and randomly drawn from a Gaussian, in the present case the diagonal correlation matrix is $\Lambda={\diag}(294.845, 34.679, 30.311, 11.612, 10.733, 9.468, 8.232, 5.307, 4.144, 2.443, 2.429, 2.218,$ $2.083, 1.986, 1.406, 1.382, 1.102, 1.001, 0.889, 0.707, 0.693, 0.684, 0.665, 0.63, 0.594, 0.591, 0.576,$ $0.574, 0.562, 0.467, 0.463, 0.455)$. The largest eigenvalues are three outliers which can be seen in the vicinity of $x=-1$, see the inset.}
 \label{figsimsaddlepoint}
\end{figure}

We compared the result~\eqref{densitySaddlePoint3} with Monte Carlo simulations in Fig.~\ref{figsimsaddlepoint}. We find a perfect agreement for a matrix size $p=32$. Due to the weaker level repulsion of real matrices compared to complex ones the agreement is in the bulk better for $\beta=1$ than for $\beta=2$ while it is worse in the tails at the edges. Indeed the strong oscillations for the complex ensembles is a direct result of the stronger level repulsion which is quadratic. Even the positions of the outliers are astoundingly good predicted and the shape of their distributions moderately approximated.

We want to emphasize that we obtain the well-known results~\cite{BJNPZ} for the uncorrelated case. Then the distribution of the Jacobi ensemble has either  square root zeros or square root singularities at the boundaries of the support depending on if the edges detach from the generic bounds $x=\pm1$ or not. This can be easily checked by our results. For the Cauchy-Lorentz ensemble we obtain a Levy-tail with an algebraic decay of $b^{-3/2}$ if the corresponding Jacobi ensemble does not detach from the lower bound $x=-1$ which was also found in \cite{BJNPZ}.

Finally, let us come back to the case of a degeneracy of $C_{\rm eff}$ and, hence, $\Lambda$. When taking $p\to lp$, $n_1\to ln_1$, $n_2\to ln_2$, and $\Lambda\to \Lambda\otimes\eins_l$, we immediately notice that the number of copies $l\in\mathbb{N}$ drops out in the saddle point equation~\eqref{saddle}. Thus the saddlepoint is independent of $l$. Also in the level densities~\eqref{densitySaddlePoint2} and \eqref{densitySaddlePoint3} the number $l$ does not appear. Hence the asymptotic result is independent. We expect that the local statistics remain unaffected, too, as it is the case for the correlated Wishart ensemble, see \cite{WirtzKieburgGuhrPRL}. For the complex case this can be readily checked via the kernel~\eqref{projection-res}.

\section{Conclusions}\label{sec:conclusion}

We considered the correlated Jacobi and Cauchy-Lorentz ensemble.
Our first result is the map of the correlated Jacobi ensemble to the correlated Cauchy-Lorentz ensemble. The eigenvalue statistics of the first is completely determined by the latter and \textit{vice versa}. This is impressive since the first ensemble exhibits a spectrum on a compact interval namely $[-1,1]$ while the second one has a Levy tail which drops off in the uncorrelated case with the algebraical behavior $x^{3/2}$. The uncorrelated Cauchy-Lorentz ensemble was tried to apply to finance~\cite{BJNPZ} because of the known heavy tail behavior in the statistics.

We also derived a supersymmetric integral for the $k$-point generating function of the correlated Jacobi ensemble and the correlated Cauchy-Lorentz ensemble via the projection formula~\cite{KaymakKieburgGuhr,Kie15} in combination with generalized Hubbard-Stratonovich transformation~\cite{IngSie,KGG09,KSG} and the superbosonization formula~\cite{KSG,LSZ,Sommers}. The resulting integral over the supermatrices looks similar to the one of the correlated Wishart ensemble~\cite{Rech12}. This representation is ideal for studying the asymptotic limit ($p\sim n_1\sim n_2\to\infty$) due to its small and fixed number of integration variables. For example we calculate a closed expression for the limit of the macroscopic level density. The approximation is already very good for moderate matrix size $p\sim 30$  and generic empirical fixed correlations which is confirmed by Monte Carlo simulations. We underline that even the outliers are predicted by the asymptotic limit.

The level densities at finite matrix dimensions were explicitly calculated for both the correlated real and complex random matrix ensembles. As in the Wishart case the real ensembles are more involved. Nevertheless we could simplify the result to a finite sum where each summand is a product of three integrals and one of the integrals can be performed exactly. Both technical properties resemble the results for the Wishart ensemble, see~\cite{Rech12}. The remaining one-fold integrals can be numerically evaluated. In the case of double degeneracy of the spectrum of the empirical covariance matrix all integrals can be analytically performed. This generic degeneracy was proposed in~\cite{WirtzKieburgGuhrPRL} where it was shown that the asymptotic spectral statistics of correlated Wishart ensembles do not differ from the case without degeneracy.  We also observed that the degeneracy has no influence on the asymptotic statistics of the correlated Jacobi ensemble and the correlated Cauchy-Lorentz ensemble. Thus we propose also to 
study the case with a doubly degenerated covariance matrix artificially introduced by taking two copies of the covariance matrix in the case of real ensembles due to its analytical advantage.

In the real case we had to restrict our explicit calculation to the level density. Higher order correlations are not analytical feasible for finite matrix dimensions at the moment due to the lack of knowledge about certain group integrals. This is not the case for the complex ensembles where we took a different approach as for the real matrices. For the correlated complex Jacobi ensemble and Cauchy-Lorentz ensemble we first derived the joint probability density of their eigenvalues. This was possible due to an Itzykson-Zuber-Harish-Chandra-like group integral derived in~\cite{GroRic,HarOrl}. The joint probability density satisfies a determinantal point process and the corresponding kernel resembles the result of the correlated Wishart ensemble~\cite{SimMou,ATLV}. Nonetheless we derived also for this kernel a supersymmetric integral which is much more suitable to study the behavior at large matrix sizes.

Our analysis can be extended into various directions: First of all, the supersymmetric expression for the $k$-point correlation function is a perfect starting point to derive also for $k>1$ closed expressions. The  distributions of the largest and the smallest eigenvalues are other important quantities which were already studied for the correlated Wishart ensemble~\cite{For06,WKG14,WirGuh}.  Also the case of double correlations, as discussed in Ref.~\cite{WaltnerWirtzGuhr2014} for Wishart ensembles, could be considered. Another generalization could be the investigation of other correlated heavy tailed ensembles instead of the Cauchy-Lorentz one. Product matrices, see~\cite{AkeIps} for a recent review, yield a new and analytical feasible approach to such heavy tailed ensembles. The combination of the approach applied in the present work with product matrices and the projection formula~\cite{KaymakKieburgGuhr,Kie15} may provide a unique and ideal tool to study the macroscopic as well as the local spectral 
statistics of Levy tailed ensembles.

\section{Acknowledgments}
T.W. acknowledges support from the German Research Council (DFG) via the Sonderforschungsbereich Transregio 12, ``Symmetries and Universality in Mesoscopic Systems''. M.K. partially acknowledges financial support from the Alexander von Humboldt-Foundation and from  the CRC 701: \textit{Spectral Structures and Topological Methods in Mathematics} of the DFG.

\appendix

\section{Supersymmetric two-Matrix Model}\label{sec:cje:twosupermatrixmodel}

 Let us derive another supersymmetric integral which consists of two supermatrices entering in a symmetric way.  This integral explicitly shows the symmetry of the correlated Jacobi ensemble under $n_1\leftrightarrow n_2$, $\Lambda\to\Lambda^{-1}$, and $x\to -x$. This symmetry is not immediate in the expression~\eqref{eq:cje:chifinal} where we have to substitute $b=(1-x)/(1+x)$ since it is given for the correlated Cauchy ensemble. However the supersymmetric integral~\eqref{eq:cje:chifinal} is certainly simpler to compute than the one we present in this section since we have to deal with only one supermatrix in Eq.~\eqref{eq:cje:chifinal}.
 
We start from Eq.~(\ref{eq:cje:reminderGeneratingFunction}). To apply the same approach as in \cite{Rech12} we have to linearize the arguments of the characteristic polynomials in the $FF^\dagger$ and $BB^\dagger$. This can be achieved by
multiplying the matrices in the determinants from the right with $(FF^\dagger + BB^\dagger)$ yielding
\begin{align}
  \label{eq:cje:linearizingratio}
\frac{\det\left((FF^\dagger - BB^\dagger)(FF^\dagger + BB^\dagger)^{-1} - \kappa_{a2}\eins_{p}\right)}{\det\left((FF^\dagger - BB^\dagger)(FF^\dagger + BB^\dagger)^{-1}-
\kappa_{b1}\eins_{p}\right)} = \left(\frac{1+\kappa_{a2}}{1+\kappa_{b1}}\right)^p\frac{\det\left(FF^\dagger\widehat{\kappa}_{a2} - BB^\dagger\right)}
{\det\left(FF^\dagger\widehat{\kappa}_{b1} - BB^\dagger\right)}.
\end{align}
with $\widehat{\kappa}=(1-\kappa)/(1+\kappa)$. In the next step we plug Eq.~\eqref{eq:cje:linearizingratio} into Eq.~(\ref{eq:cje:reminderGeneratingFunction}) and express the determinants as a Gaussian integral over a rectangular supermatrix
\begin{eqnarray}
 \label{eq:gHS:Abeta1}A&=&\left[z_{ja}\,z_{ja}^*\,\zeta_{jb}\,\zeta_{jb}^*\right],\ \beta=1,\\
 \label{eq:gHS:Abeta2}A&=&\left[z_{ja}\,\zeta_{jb}\right],\ \beta=2,
\end{eqnarray}
 of dimension $p\times(\gamma k|\gamma k)$, \textit{i.e.}
\begin{align}
  \label{eq:cje:expressingdetsasgaussians}
\prod_{a=1}^k\frac{\det\left(FF^\dagger\widehat{\kappa}_{a2} - BB^\dagger\right)}{\det\left(FF^\dagger\widehat{\kappa}_{a1} - BB^\dagger\right)} &=
\int d[A]\exp\left(\imath \tr FF^\dagger A\jmat A^\dagger + \imath\tr BB^\dagger AA^\dagger\right).
\end{align}
Here we assume for simplicity that the imaginary parts of $\widehat{\kappa}_1={\rm diag}(\widehat{\kappa}_{11},\ldots,\widetilde{\kappa}_{k1})$ are on the complex upper half-plane. The source matrix is
 \begin{align}
\label{eq:cje:introducingJmat}
  \jmat ={{\diag}}\left(\widehat{\kappa}_{11},\ldots,\widetilde{\kappa}_{k1},\widehat{\kappa}_{12},\ldots,\widetilde{\kappa}_{k2}\right)\otimes\eins_\gamma.
\end{align}
We substitute the integral (\ref{eq:cje:expressingdetsasgaussians}) into the generating function (\ref{eq:cje:reminderGeneratingFunction}) and exchange the $F$ and $B$  with the $A$ integral.
The  resulting $F$ and $B$ integrals are Gaussian and yield
\begin{eqnarray}
  \label{eq:cje:characteristicfunctionWishart}
&&\int d[F]\int d[B]P(F|\eins_p)P(B|\Lambda) \exp\left(\imath \tr FF^\dagger A\jmat A^\dagger + \imath\tr BB^\dagger AA^\dagger\right)\\
& =& {\det}^{-n_1/\gamma} \left(\eins_{p}-\imath A\jmat A^\dagger\right){\det}^{-n_2/\gamma} \left(\eins_{p}-\imath \Lambda A A^\dagger\right).\nonumber
\end{eqnarray}
Then the generating function becomes
\begin{align}
    Z_{p,\beta}^{k|k}(\kappa) \propto{\sdet}^{-p}\left(\eins_{k}+\kappa\right)\int d[A]
    {\det}^{-n_1/\gamma}\left(\eins_{p}-\imath
      A\jmat A^\dagger\right){\det}^{-n_2/\gamma}\left(\eins_{p}-\imath
      \Lambda AA^\dagger\right). \label{eq:cje:noABgeneratingfunction}
\end{align}
The normalization constant is independent of $\kappa$ and $\Lambda$.
The next step is known as the duality between ordinary and superspace. Due to the invariance of the integrand in Eq.~(\ref{eq:cje:noABgeneratingfunction}) under $A\to UA$ for an arbitrary $U\in{\rm O}(p)$ for $\beta=1$ and  $U\in{\rm U}(p)$ for $\beta=2$, the integrand only depends on the invariants $\tr \left(AA^\dagger\right)^m$ for $m\in \mathbb N$. These invariants are equal to
the superinvariants $\tr \left(AA^\dagger\right)^m=\str \left( A^\dagger A\right)^m$, see \cite{Efebook,Guhrbook,IngSie,KGG09}. Employing this duality in the generating function~(\ref{eq:cje:noABgeneratingfunction}), we arrive at
\begin{align}
    Z_{p,\beta}^{k}(\kappa)\propto\sdet^{-p}\left(\eins_{k}+\kappa\right)\int d[A] {\sdet}^{-n_1/\gamma_1}\left(\eins_{{\gamma} k|\gamma k}-\imath A^\dagger A \jmat\right)
    {\sdet}^{-n_2/\gamma_1}\left(\eins_{\gamma k|\gamma k}-\imath A^\dagger\Lambda A\right).
\label{eq:cje:finalFFexpression}
\end{align}

The main difference of Eq.~\eqref{eq:cje:finalFFexpression} to most models discussed in the literature so far is that this one includes two different products of $A$ and $A^\dagger$. Namely, $A^\dagger A $ which arises naturally if invariant matrix models are considered and $A^\dagger \Lambda_{\text{eff}}A$ appearing due to a non-trivial correlation structure. We cannot replace both products by one supermatrix, but we can apply the generalized Hubbard-Stratonovich transformation~\cite{IngSie,KGG09,KSG} independently for both products. It yields
the following supersymmetric two-matrix model
\begin{align}
  \begin{split}
      Z_{p,\beta}^{k}(\kappa) &\propto\sdet^{-p}\left(\eins_{k}+\kappa\right)\int
       d[\sigma]d[\varrho]
       I_{n_2}(\varrho)I_{n_1}(\sigma)\exp\left(-\str\varrho-\str\sigma\right)\\&\times
       \sdet^{-1/\gamma}\left(\eins_p\otimes\sigma -
         \Lambda^{-1}\otimes\varrho\jmat \right),
    \end{split}
    \label{eq:cje:twosupermatrixmodel}
\end{align}
where the function $I_{n_i}(\varrho)$, $i=1,2$, is the supersymmetric Ingham-Siegel integral~\eqref{eq:gHS:InghamSiegel}. The $(\gamma k|\gamma k)\times(\gamma k|\gamma k)$ dimensional supermatrices $\rho$ and $\sigma$ have the same symmetries as the supermatrix in the third equality of Eq.~\eqref{eq:cje:chifinal}.

We can completely symmetrize the integral in $\rho$ and $\sigma$ by going back to the sources $\widehat{\kappa}\to\kappa$ and the empirical matrix $\Lambda \to C_F^{-1/2}C_BC_F^{-1/2}$. Then we have the final result
\begin{align}
  \begin{split}
      Z_{p,\beta}^{k}(\kappa) &\propto\int
       d[\sigma]d[\varrho]
       I_{n_1}(\varrho)I_{n_2}(\sigma)\exp\left(-\str\varrho-\str\sigma\right)\\
       &\times
       \sdet^{-1/\gamma}\left(C_B\otimes\sigma[(1+\kappa)\otimes\eins_\gamma] -
         C_F\otimes\varrho[(1-\kappa)\otimes\eins_\gamma] \right).
    \end{split}
    \label{eq:cje:twosupermatrixmodel-fin}
\end{align}
This expression is completely invariant under the original symmetry $n_1\leftrightarrow n_2$, $\Lambda\to\Lambda^{-1}$, and $\kappa\to -\kappa$ because the symmetry is achieved by the change $\rho\leftrightarrow\sigma$.

We again underline that the supermatrix model~(\ref{eq:cje:twosupermatrixmodel}) can be in principal computed by expanding the integrand in the Grassmann variables and performing the remaining integrals. However we have now two supermatrices such that this calculation can be a highly non-trivial task. This is the reason why we use more advanced techniques which involve the relation to the correlated Cauchy-Lorentz ensemble, see section~\ref{sec:JCL}. The supersymmetry result~\eqref{eq:cje:chifinal} can be obtained from Eq.~\eqref{eq:cje:twosupermatrixmodel} by rescaling $\sigma\to\rho\sigma$ and then integrating over $\rho$ which yields the superdeterminant $\sdet^{-\mu/\gamma}(\sigma+\eins_{\gamma k|\gamma k})$.

\section{Regularizations of the Integrals in Section~\ref{sec:cje:realcase}}\label{app:expl}

The numerical evaluation of the integrals~\eqref{int2} and \eqref{int3} suffer by the non-integrable singularities of order $3/2$ at the boundaries, in particular they are of the two forms
\begin{align}
 J_1=\int\limits_{b\Lambda_{j+1}^{-1}}^{b\Lambda_j^{-1}}d{r}\frac{f(r)}{\left|b\Lambda_j^{-1}-r\right|^{3/2}},\quad J_2=\int\limits_{b\Lambda_j^{-1}}^{b\Lambda_{j-1}^{-1}}d{r}\frac{f(r)}{\left|b\Lambda_j^{-1}-r\right|^{3/2}}\label{eq:cje:Rdensity:example2.I}
\end{align}
for certain real valued functions $f(r)$ without singularities in the interval $[b\Lambda_{j+1}^{-1},b\Lambda_j^{-1}]$. As already said the integrals are taken via Cauchy's principal value because of the original imaginary increment $\imath\varepsilon$. Thus we can effectively regularize the integral as follows~\footnote{We thank Petr Braun for showing us this technical trick.}
\begin{align}
\begin{split}
J_1&= \lim_{\varepsilon\to0}{\rm Re}\int\limits_{b\Lambda_{j+1}^{-1}}^{b\Lambda_j^{-1}+\varepsilon}d{r}\frac{f(r)}{\left((b+\imath\varepsilon)\Lambda_j^{-1}-r\right)^{3/2}}\\
&=\lim_{\varepsilon\to0}{\rm Re}\int\limits_{b\Lambda_{j+1}^{-1}}^{b\Lambda_j^{-1}+\varepsilon}d{r}\frac{f(r)-f(b\Lambda_{j}^{-1})}{\left((b+\imath\varepsilon)\Lambda_j^{-1}-r\right)^{3/2}}+\lim_{\varepsilon\to0}{\rm Re}\left[\frac{2f(b\Lambda_{j}^{-1})}{\sqrt{(b+\imath\varepsilon)\Lambda_j^{-1}-r}}\right]_{r=b\Lambda_{j+1}^{-1}}^{b\Lambda_j^{-1}+\varepsilon}\\
&=\int\limits_{b\Lambda_{j+1}^{-1}}^{b\Lambda_j^{-1}}d{r}\frac{f(r)-f(b\Lambda_j^{-1})}{\left|b\Lambda_j^{-1}-r\right|^{3/2}}-\frac{2f(b\Lambda_{j}^{-1})}{\sqrt{|b\Lambda_j^{-1}-b\Lambda_{j+1}^{-1}|}}
\end{split}\label{eq:cje:Rdensity:example2.II}
\end{align}
and similar for the other integral (then the imaginary part is needed)
\begin{align}
\begin{split}
J_2&= \int\limits_{b\Lambda_j^{-1}}^{b\Lambda_{j-1}^{-1}}d{r}\frac{f(r)-f(b\Lambda_j^{-1})}{\left|b\Lambda_j^{-1}-r\right|^{3/2}}-\frac{2f(b\Lambda_{j}^{-1})}{\sqrt{|b\Lambda_j^{-1}-b\Lambda_{j-1}^{-1}|}}.
\end{split}\label{eq:cje:Rdensity:example2.II.b}
\end{align}
The minus sign in Eq.~\eqref{eq:cje:Rdensity:example2.II.b} in front of the second term results from taking the imaginary part despite it is evaluated at the upper boundary. The cut-off of the intervals can be also chosen independently of $\varepsilon$. The reason is that the other boundary term of the integration by parts vanishes due to taking the real or imaginary part, respectively.

We define the following two one-fold integrals
\begin{equation}\label{int-g1}
g_{a,c,l}^0(b;\Lambda)=\int_{V_l}\frac{r^{a/2}dr}{(1+r)^{c/2}\sqrt{|\det(b\Lambda^{-1}-r\eins_p)|}}
\end{equation}
and
\begin{eqnarray}
&&g_{a,c,l}^1(b;\Lambda;\Lambda_i)\label{int-g2}\\
&=&\left\{\begin{array}{cl}
\displaystyle {\rm sign}(p-l-i)\int_{V_l}\frac{r^{a/2}dr}{(1+r)^{c/2}|b\Lambda_i^{-1}-r|\sqrt{|\det(b\Lambda^{-1}-r\eins_p)|}}, & \\
& \hspace*{-3.5cm} l\neq p-i,p-i+1,\\
&\\
\displaystyle-\frac{2(b\Lambda_{i}^{-1})^{a/2}}{(1+b\Lambda_{i})^{c/2}\sqrt{|b\Lambda_i^{-1}-b\Lambda_{i+1}^{-1}|}\sqrt{|\det(b\Lambda_{\neq i}^{-1}-b\Lambda_{i}^{-1}\eins_{p-1})|}}+\int_{V_{p-i}}\frac{dr}{|b\Lambda_i-r|^{3/2}}& \\
\displaystyle\times\left[\frac{r^{a/2}}{(1+r)^{c/2}\sqrt{|\det(b\Lambda_{\neq i}^{-1}-r\eins_{p-1})|}}-\frac{(b\Lambda_i^{-1})^{a/2}}{(1+b\Lambda_i^{-1})^{c/2}\sqrt{|\det(b\Lambda_{\neq i}^{-1}-b\Lambda_i^{-1}\eins_{p-1})|}}\right],&\\
&\hspace*{-3.5cm} l=p-i,\\
&\\
\displaystyle\frac{2(b\Lambda_{i}^{-1})^{a/2}}{(1+b\Lambda_{i})^{c/2}\sqrt{|b\Lambda_i^{-1}-b\Lambda_{i-1}^{-1}|}\sqrt{|\det(b\Lambda_{\neq i}^{-1}-b\Lambda_{i}^{-1}\eins_{p-1})|}}-\int_{V_{p-i+1}}\frac{dr}{|b\Lambda_i-r|^{3/2}}& \\
\displaystyle\times\left[\frac{r^{a/2}}{(1+r)^{c/2}\sqrt{|\det(b\Lambda_{\neq i}^{-1}-r\eins_{p-1})|}}-\frac{(b\Lambda_i^{-1})^{a/2}}{(1+b\Lambda_i^{-1})^{c/2}\sqrt{|\det(b\Lambda_{\neq i}^{-1}-b\Lambda_i^{-1}\eins_{p-1})|}}\right],&\\
&\hspace*{-3.5cm} l=p-i+1.
\end{array}\right.\nonumber
\end{eqnarray}
Then we can combine the discussion about the splitting of the integral over $\mathbb{R}_+^2$ into disjoint sets and the regularization of the $3/2$-singularities. Therefore we explicitly have for the imaginary parts of the integrals~\eqref{int1}, \eqref{int2}, and \eqref{int3}
\begin{eqnarray}
\frac{1}{\pi}\lim_{\varepsilon\to0}{\rm Im}\,N_{a,c,d_1,d_2}^0(b+\imath\varepsilon;\Lambda)&=&\frac{1}{\pi}\sum_{\substack{0\leq l_1,l_2\leq p \\ l_1+l_2\in2\mathbb{N}_0+1}}(-1)^{(l_1+l_2+1)/2}{\rm sign}(l_1-l_2)\\
&&\times\det\left[\begin{array}{cc} g_{a+d_1+1,c+d_1+1,l_1}^0(b;\Lambda) & g_{a+d_1-1,c+d_1+1,l_1}^0(b;\Lambda) \\ g_{a+d_2+1,c+d_2+1,l_2}^0(b;\Lambda) & g_{a+d_2-1,c+d_2+1,l_2}^0(b;\Lambda) \end{array}\right],\nonumber\\
\frac{1}{\pi}\lim_{\varepsilon\to0}{\rm Im}\,N_{a,c,d}^1(\kappa;\Lambda;\Lambda_j)&=&\frac{1}{\pi}\sum_{\substack{0\leq l_1,l_2\leq p \\ l_1+l_2\in2\mathbb{N}_0+1}}(-1)^{(l_1+l_2+1)/2}{\rm sign}(l_1-l_2)\\
&&\times\det\left[\begin{array}{cc} g_{a+2,c,l_1}^1(b;\Lambda;\Lambda_j) & g_{a,c,l_1}^1(b;\Lambda;\Lambda_j) \\ g_{a+d+1,c+d+1,l_2}^0(b;\Lambda) & g_{a+d-1,c+d+1,l_2}^0(b;\Lambda) \end{array}\right],\nonumber\\
\frac{1}{\pi}\lim_{\varepsilon\to0}{\rm Im}\,N_{a,c}^2(\kappa;\Lambda;\Lambda_i,\Lambda_j)&=&\frac{1}{\pi}\sum_{\substack{0\leq l_1,l_2\leq p \\ l_1+l_2\in2\mathbb{N}_0+1}}(-1)^{(l_1+l_2+1)/2}{\rm sign}(l_1-l_2)\\
&&\times\det\left[\begin{array}{cc} g_{a+2,c,l_1}^1(b;\Lambda;\Lambda_i) & g_{a,c,l_1}^1(b;\Lambda;\Lambda_i) \\ g_{a+2,c,l_2}^1(b;\Lambda;\Lambda_j) & g_{a,c,l_2}^1(b;\Lambda;\Lambda_j) \end{array}\right].\nonumber
\end{eqnarray}
These results can be combined with Eq.~\eqref{eq:cje:Rdensity:derZ} to find the level density $S'_1(b)$ of the correlated Lorentz ensemble, see Eq.~\eqref{CL-dens-result}.

\end{document}